\documentclass[preprint,12pt,authoryear]{elsarticle}

\usepackage{amssymb}

\usepackage{amsthm,amssymb,amsbsy,latexsym,amsmath}
\usepackage{amsfonts}
\usepackage{amsmath}
\usepackage{graphicx}
\usepackage{float,multirow} 
\usepackage{pgf,tikz}

\theoremstyle{plain}
\newtheorem{theorem}{Theorem}[section]
\newtheorem{lemma}[theorem]{Lemma}
\newtheorem{corollary}[theorem]{Corollary}
\newtheorem{example}[theorem]{Example}
\newtheorem{prop}[theorem]{Proposition}
\theoremstyle{definition}

\journal{Discrete Mathematics}

\begin{document}

\begin{frontmatter}



\title{Formulas for the Number of Weak Homomorphisms from Paths to Rectangular Grid Graphs}

\author[mymainaddress]{Penying Rochanakul}
\ead{penying.rochanakul@cmu.ac.th}

\author[mymainaddress]{Hatairat Yingtaweesittikul}
\ead{hatairat.y@cmu.ac.th}

\author[mymainaddress]{Sayan Panma\corref{mycorrespondingauthor}}
\cortext[mycorrespondingauthor]{Corresponding author}
\ead{sayan.panma@cmu.ac.th}

\address[mymainaddress]{Department of Mathematics, Faculty of Science,
	Chiang Mai University, Chiang Mai, Thailand }

\begin{abstract}
A \emph{weak homomorphism} from a graph $G$ to a graph $H$ is a mapping $f:V(G)\to V(H)$, where either $f(x) = f(y)$ or $\{f(x), f(y)\}\in E(H)$ holds for all $\{x, y\} \in E(G)$. A  \emph{rectangular grid graph} is formed by taking the Cartesian product of two paths. In this paper, we present a formula for calculating the number of weak homomorphisms from paths to rectangular grid graphs.
\end{abstract}

%


\end{frontmatter}


\section{Introduction}
\label{1}

\hspace{0.1cm} In mathematics, the image refers to the set of values obtained by applying a mapping to all elements within the domain. Within this image, certain structural properties of the domain are retained. A mapping that maintains such a structure, which is of particular interest for our study, is commonly referred to as a homomorphism. In the context of graphs, a homomorphism is defined as follows.

 Consider the graphs $G$ and $H$. A mapping, denoted as $f : V(G) \to V(H),$ is a  \emph{homomorphism} from $G$ to $H$ if $\{f(x), f(y)\} \in E(H)$ for all $\{x, y\} \in E(G)$, meaning that $f$ \emph{preserves} the edges.
The set of homomorphisms from $G$ to $H$ is denoted as $\mathrm{Hom}(G, H)$.
Let $P_{n}$ represent a path of order $n$ with vertex set $V(P_{n}) = \{0, 1, ..., n-1\}$ and edge set $E(P_{n}) = \{\{i, i+1\} | i = 0, 1, ..., n-2\}$.
Similarly, $C_{n}$ stands for a cycle of order $n$ ($n \geqslant 3$) with vertices $V(C_{n}) = \{0, 1, ..., n-1\}$ and edges $E(C_{n}) = \{\{i, i+1\} | i = 0, 1, ..., n-1\}$, where the addition is performed modulo $n$.
For a deeper understanding of graphs and algebraic graphs, we direct readers to references \cite{4} and \cite{61}. 

The expression for determining the number of homomorphisms from $P_n$ to itself, known as End($P_n$), was introduced by Arworn in 2009 \cite{1}. Arworn transformed the problem by equating it to enumerating the shortest paths originating from point $(0,0)$ and reaching any point $(i,j)$ within an $r$-ladder square lattice, ultimately deriving a succinct formula.

In a broader context, a homomorphism from a graph $G$ to itself is termed an \emph{endomorphism} on $G$. It is evident that the set of endomorphisms on $G$ constitutes a monoid, wherein the composition of mappings serves as the defining operation.

When considering a mapping $f: V(G) \to V(H)$, the concept of $f$ \emph{contracting} an edge $\{x, y\}$ denotes that both vertices $x$ and $y$ are mapped to the same vertex in $V(H)$, i.e., $f(x) = f(y)$.
The central concept is that homomorphisms must preserve edges.
If we also have the option to contract edges, then this achievement can be realized using regular homomorphisms when our graphs contain a loop at every vertex.

A mapping $f : V(G) \to V(H)$ is termed a \emph{weak homomorphism} from a graph $G$ to a graph $H$ (also referred to as an \emph{egamorphism}) if $f$ contracts or preserves the edges, that is, $f(x) = f(y)$ or $\{f(x), f(y)\} \in E(H)$ whenever $\{x, y\} \in E(G)$.
A weak homomorphism from $G$ to itself is referred to as a \emph{weak endomorphism} on $G$.
We denote the set of weak homomorphisms from $G$ to $H$ as WHom($G, H$) and the set of weak endomorphisms on $G$ as WEnd($G$).
It is evident that WEnd($G$) constitutes a monoid under the composition of mappings. The composition of (weak) homomorphisms also forms a (weak) homomorphism. Consequently, this results in a preorder on graphs and defines a category \cite{4}.

In 2010, Sirisathianwatthana and Pipattanajinda \cite{5} established the count of weak homomorphisms of cycles as WHom($C_m, C_n$), expressed in terms of the collection of WHom$^i_j$($P_{m-1}, C_n$), where WHom$^i_j$($P_{m-1}, C_n$) represents a set of weak homomorphisms from $P_{m-1}$ to $C_n$, with the conditions that $f(0)=i$ and $f(m-1)=j$. In 2018, Knauer and Pipattanajinda \cite{7} introduced the count of weak endomorphisms on paths, denoted as WEnd($P_n$), by relating it to the quantities of shortest paths from the origin point $(0,0,0)$ to any arbitrary point $(i,j,k)$ within the three-dimensional square lattice, as well as within the $r$-ladder three-dimensional square lattice. Moreover, they provided formulas for the count of shortest paths from the point $(0,0,0)$ to any point $(i,j,k)$, as shown in Proposition \ref{prop1}. Figures \ref{Figure1} and \ref{Figure1A} depict the cubic lattice and the $2$-ladder cubic lattice when $i=6$, $j=4$, and $k=4$, respectively.

\begin{prop}[\cite{7}]\label{prop1} The numbers $M(i, j, k)$ and $M_r(i, j, k)$ of shortest
	paths from the point $(0, 0, 0)$ to any point $(i, j, k)$ in the  cubic lattice and in the $r$-ladder  cubic lattice are
	$$M(i,j,k) =  \dbinom{i+j+k}{i,j,k}$$
	and
	$$M_r(i,j,k) =  \left[ \dbinom{i+j+k}{i,j,k} - \dbinom{i+j+k}{j-r-1,i+r+1,k}  \right], $$
	respectively.
\end{prop}

\begin{figure}[h!]
	\begin{minipage}[b]{0.5\textwidth}	
		\begin{tikzpicture}[scale=0.4]
			\def\x{6}
			\def\y{4}
			\def\z{4}
			\foreach \i in {0,1,...,\y} 
			\draw[line width=0.3mm,gray] (0, \i) -- ({1.2*\x},{\i+0.2*\x});
			\foreach \i in {1,...,\z} 
			\draw[line width=0.3mm,gray] ({-1*\i},{0.3*\i+\y}) -- ({-1*\i+1.2*\x},{0.3*\i+0.2*\x+\y});
			\foreach \j in {0,1,...,\x}
			\draw[line width=0.6mm,black] ({1.2*\j},{0.2*\j}) -- ({1.2*\j},{0.2*\j+\y});
			\foreach \j in {1,...,\z}
			\draw[line width=0.6mm,black] ({-1*\j},{0.3*\j}) -- ({-1*\j},{0.3*\j+\y});
			\foreach \k in {0,...,\y}
			\draw[dashed,magenta] (0,\k) -- ({-1*\z},{0.3*\z+\k});
			\foreach \k in {1,...,\x}
			\draw[dashed,magenta] ({1.2*\k},{\y+0.2*\k}) -- ({1.2*\k-\z},{\y+0.2*\k+0.3*\z});
			\foreach \i in {0,...,\y} {
				\foreach \j in {0,1,...,\x}
				\draw[fill=blue] ({1.2*\j},{\i+0.2*\j}) circle(1mm);
				\foreach \k in {1,...,\z}
				\draw[fill=magenta] ({-1*\k},{\i+0.3*\k}) circle(1mm);}
			\foreach \i in {1,...,\x} {
				\foreach \j in {1,...,\z}
				\draw[fill=green] ({1.2*\i-\j},{\y+0.2*\i+0.3*\j}) circle(1mm);
			}
			\draw (0,-0.7) node[inner sep=1pt,rectangle,fill=white] {$(0,0,0)$};
			\draw (1.2*\x,{-1+0.2*\x}) node[inner sep=1pt,rectangle,fill=white] {$(6,0,0)$};
			\draw ({-1*\z},{-1+0.3*\z}) node[inner sep=1pt,rectangle,fill=white] {$(0,0,4)$};
			\draw ({1.2*\x-\z},{0.7+\y+0.3*\z+0.2*\x}) node[inner sep=1pt,rectangle,fill=white] {$(6,4,4)$};
			\draw ({1.2*\x},{1+\y+0.2*\x}) node[inner sep=1pt,rectangle,fill=white] {$(6,4,0)$};
			\draw ({-1*\z},{1+\y+0.3*\z}) node[inner sep=1pt,rectangle,fill=white] {$(0,4,4)$};
		\end{tikzpicture}	\caption{\label{Figure1} Cubic lattice }
	\end{minipage}
	\begin{minipage}[b]{0.5\textwidth}
		\begin{tikzpicture}[scale=0.4]
			\def\x{6}
			\def\y{4}
			\def\z{4}
			\def\r{2}
			\def\a{2}
			\def\b{4}
			\foreach \i in {0,...,\y} 
			\draw[line width=0.3mm,gray] (0, \i) -- ({1.2*\x},{\i+0.2*\x});
			\foreach \i in {0,...,\a} 
			\draw[line width=0.5mm,white] (0, {\r+\i}) -- ({1.2*\i},{\r+1.2*\i});
			\foreach \i in {1,...,\z} 
			\draw[line width=0.3mm,gray] ({-1*\i+1.2*\a},{0.2*\a+0.3*\i+\y}) -- ({-1*\i+1.2*\x},{0.3*\i+0.2*\x+\y});
			\foreach \j in {1,...,\a} {
				\foreach \k in {1,...,\z} 
				\draw[line width=0.3mm,gray] ({1.2*\j-\k},{\j+\r-1+0.2*\j+0.3*\k}) -- ({1.2*\j-\k-1.2},{\j+\r-1.2+0.2*\j+0.3*\k}) ;}
			\foreach \j in {1,...,\b}
			\draw[line width=0.6mm,black] ({1.2*\j+1.2*\a},{0.2*\j+0.2*\a}) -- ({1.2*\j+1.2*\a},{0.2*\j+\y+0.2*\a});
			\foreach \j in {0,1,...,\a}
			\draw[line width=0.6mm,black] ({1.2*\j},{0.2*\j}) -- ({1.2*\j},{0.2*\j+\r+\j});
			\foreach \j in {1,...,\z}
			\draw[line width=0.6mm,black] ({-1*\j},{0.3*\j}) -- ({-1*\j},{0.3*\j+\r});
			\foreach \j in {1,...,\a} {
				\foreach \k in {1,...,\z} 
				\draw[line width=0.6mm,black] ({1.2*\j-\k},{\j+\r-1+0.2*\j+0.3*\k}) -- ({1.2*\j-\k},{\j+\r+0.2*\j+0.3*\k}) ;}
			\foreach \k in {0,...,\r}
			\draw[dashed,magenta] (0,\k) -- ({-1*\z},{0.3*\z+\k});
			\foreach \k in {1,...,\b}
			\draw[dashed,magenta] ({1.2*\k+1.2*\a},{\y+0.2*\k+0.2*\a}) -- ({1.2*\k-\z+1.2*\a},{\y+0.2*\k+0.3*\z+0.2*\a});
			\foreach \j in {1,...,\a}
			{\draw[dashed,magenta] ({1.2*\j},{1.2*\j+\r}) -- ({-\z+1.2*\j},{0.3*\z+1.2*\j+\r}) ;
				\draw[dashed,magenta] ({1.2*\j},{1.2*\j+\r-1})-- ({-1*\z+1.2*\j},{0.3*\z+1.2*\j+\r-1}) ;}
			\foreach \j in {0,...,\r} {
				\foreach \i in {0,...,\a}
				\draw[fill=blue] ({1.2*\i},{\j+0.2*\i}) circle(1mm);}
			\foreach \j in {1,...,\a} {
				\foreach \i in {1,...,\j}
				\draw[fill=blue] ({1.2*\j},{\i+\r+0.2*\j}) circle(1mm);}
			\foreach \j in {0,...,\y} {
				\foreach \i in {0,1,...,\b}
				\draw[fill=blue] ({1.2*\i+1.2*\a},{\j+0.2*\i+0.2*\a}) circle(1mm);}
			\foreach \j in {0,...,\r} {
				\foreach \k in {1,...,\z}
				\draw[fill=magenta] ({-1*\k},{\j+0.3*\k}) circle(1mm);}
			\foreach \j in {1,...,\a} {
				\foreach \k in {1,...,\z} {
					\draw[fill=yellow] ({1.2*\j-\k},{\j+\r-1+0.2*\j+0.3*\k}) circle(1mm);
					\draw[fill=yellow] ({1.2*\j-\k},{\j+\r+0.2*\j+0.3*\k}) circle(1mm);}}
			\foreach \i in {1,...,\b} {
				\foreach \j in {1,...,\z}
				\draw[fill=green] ({1.2*\i-\j+1.2*\a},{\y+0.2*\i+0.3*\j+0.2*\a}) circle(1mm);
			}
			\draw (0,-0.7) node[inner sep=1pt,rectangle,fill=white] {$(0,0,0)$};
			\draw ({1.2*\x},{-1+0.2*\x}) node[inner sep=1pt,rectangle,fill=white] {$(6,0,0)$};
			\draw ({-1*\z},{-1+0.3*\z}) node[inner sep=1pt,rectangle,fill=white] {$(0,0,4)$};
			\draw ({1.2*\x-\z},{0.7+\y+0.3*\z+0.2*\x}) node[inner sep=1pt,rectangle,fill=white] {$(6,4,4)$};
			\draw ({1.2*\x},{1+\y+0.2*\x}) node[inner sep=1pt,rectangle,fill=white] {$(6,4,0)$};
			\draw ({-1*\z-0.8},{0.7+\r+0.3*\z}) node[inner sep=1pt,rectangle,fill=white] {$(0,2,4)$};
			
		\end{tikzpicture}
		\caption{\label{Figure1A} $2$-ladder cubic lattice}
	\end{minipage}
\end{figure}

Recently, in 2022, Promsri et al. \cite{8} introduced the number of  weak homomorphisms of paths WHom($P_m, P_n$), by associating it with the order of the following three sets:  $A^{i}_{m-1, n} = \{ f \in $ WHom$(P_{m-1}, P_{n}) | f(0) = i \}$, $B^{i}_{m-1, n} = \{ f \in $ WHom$(P_{m-1}, P_{n}) | f(0) = i \text{ and } f(m-2) = 0 \}$, and  $C^{i}_{m-1, n} = \{ f \in $ WHom$(P_{m-1}, P_{n}) | f(0) = i \text{ and } f(m-2) = n-1 \}$, where $i$ ranges from $0$ to $n-1$. 	

For any two graphs $G_1$ and $G_2$, the \emph{Cartesian product} of $G_1$ and $G_2$ is the graph $G_1 \square G_2$ with vertices $V(G_1 \square G_2) = V(G_1) \times V(G_2)$, and in which $\{(a,u),(b,v)\}$ forms an edge if either $a=b$ and $\{u,v\}\in E(G_2)$, or $\{a,b\}\in E(G_1)$ and $u=v$. A \emph{rectangular grid graph} $P_n\square P_k$ represents the Cartesian product of $P_n$ and $P_k$.

We observe that a mapping $f:V(P_m)\to V(G_1\square G_2)$ is a homomorphism if and only if the sequence $f(0), f(1), ..., f(m-1)$ forms a walk in $G_1\square G_2$. Consequently, a one-to-one correspondence emerges between the set of homomorphisms $f:P_m\to G_1\square G_2$ and the set of walks consisting of $m$ vertices within $G_1\square G_2$. Similarly, we can establish a one-to-one correspondence between the set WHom($P_m, G_1\square G_2$) and the collection of partial walks with $m$ vertices in $G_1\square G_2$. Here, the \emph{partial walk} is a sequence obtained by concatenating $q$ walks, namely $W_1, W_2, ..., W_q$, for some $q \in \mathbb{N}$, and the ending vertex of $W_i$ is the same as the starting vertex of $W_{i+1}$ for all $i=1, 2, ..., q-1$.

In 2023, Yingtaweesittikul et al. \cite{9} introduced a formula to determine the count of homomorphisms from $P_m$ to $P_n\square P_k$, relating it to the order of the set of weak homomorphisms $f$ from $P_m$ to $P_n$ with $f(0)=j$, denoted as Hom$^j$($P_m, P_n$). This formula gives the solution to the problem concerning the number of walks of order in the rectangular grid graphs $P_n\square P_k$. Moreover, they provided formulas for Hom$^j$($P_m, P_n$), as shown in Theorem \ref{theo1}.

\begin{theorem}\cite{9}\label{theo1}
	Let $m,n$ be positive integers and $j$ a non-negative integer. Let $\mathcal{L}=\max\{0,\lceil\frac{m-j-1}{2}\rceil\}$ and $\mathcal{U}=\min\{m-1,\lfloor\frac{m+n-j-2}{2}\rfloor\}$. Then
	\begin{equation}\label{theo1equa}
		|\mathrm{Hom}^{j}(P_m,P_n)| = \sum_{i=\mathcal{L}}^{\mathcal{U}} \sum_{|t|\leq\lfloor \frac{m+n}{n}\rfloor}\left(\binom{m-1}{i-t(n+1)}-\binom{m-1}{i+j-t(n+1)+1}\right).
	\end{equation}		
\end{theorem}

If we let $m\leq n$, let $j < n$, and reduce all the zero term, we can obtain the following corollary.

\begin{corollary}\label{coro1}
	Let $m,n$ be positive integers and $j$ a non-negative integer such that $m\leq n$ and $j < n$.  Then
	\begin{equation}
		\begin{aligned}
			|\mathrm{Hom}^{j}(P_m,P_n)| =& \sum_{t=max\left\{0,\left\lceil\frac{j-(n-m)}{2}\right\rceil\right\}}^{\left\lceil\frac{m+j}{2}\right\rceil-1}\dbinom{m-1}{t} - \sum_{t=0}^{\left\lfloor\frac{j-(n-m)}{2}\right\rfloor-1}\dbinom{m-1}{t}\\&
			-\sum_{t=0}^{\left\lfloor\frac{m-j-1}{2}\right\rfloor-1}\dbinom{m-1}{t}.\nonumber
		\end{aligned}
	\end{equation}		
\end{corollary}
\begin{proof}
	Since $m\leq n$, $\lfloor \frac{m+n}{n}\rfloor\leq 2$. Thus, $t\in\{-2,-1,0,1,2\}$ and Equation (\ref{theo1equa}) can be reduced to 
	\begin{equation*}
		\begin{aligned}
			|\mathrm{Hom}^{j}(P_m,P_n)| =& \sum_{i=\mathcal{L}}^{\mathcal{U}} \left(\binom{m-1}{i+2n+2}-\binom{m-1}{i+j+2n+3}
			+\binom{m-1}{i+n+1}\right.\\&\hskip0.5cm \left.-\binom{m-1}{i+j+n+2}+\binom{m-1}{i}-\binom{m-1}{i+j+1} +\binom{m-1}{i-n-1}\right.\\&\hskip0.5cm \left.-\binom{m-1}{i+j-n}+\binom{m-1}{i-2n-2}-\binom{m-1}{i+j-2n-1}\right).	\end{aligned}
	\end{equation*}
	Since $\binom{m-1}{i+2n+2}$, $\binom{m-1}{i+j+2n+3}$, $
	\binom{m-1}{i+n+1}$, $\binom{m-1}{i+j+n+2}$, $ \binom{m-1}{i-n-1}$, $\binom{m-1}{i-2n-2}$, and $\binom{m-1}{i+j-2n-1}$ are all zeros, we have
	
	\begin{equation*}
		\begin{aligned}
			|\mathrm{Hom}^{j}(P_m,P_n)| &=\sum_{i=\mathcal{L}}^{\mathcal{U}} \left(\binom{m-1}{i}-\binom{m-1}{i+j+1}-\binom{m-1}{i+j-n}\right)\\
			&=\sum_{i=\mathcal{L}}^{\mathcal{U}} \binom{m-1}{i}-\sum_{i=\mathcal{L}}^{\mathcal{U}}\binom{m-1}{i+j+1}-\sum_{i=\mathcal{L}}^{\mathcal{U}}\binom{m-1}{i+j-n}\\
			&=\sum_{i=\mathcal{L}}^{\mathcal{U}} \binom{m-1}{(m-1)-i}-\sum_{i=\mathcal{L}}^{\mathcal{U}}\binom{m-1}{(m-1)-i-j-1}\\&\hskip0.5cm -\sum_{i=\mathcal{L}}^{\mathcal{U}}\binom{m-1}{i+j-n}\\
			&=\sum_{t=(m-1)-\mathcal{U}}^{(m-1)-\mathcal{L}} \binom{m-1}{t}-\sum_{t=(m-1)-j-1-\mathcal{U}}^{(m-1)-j-1-\mathcal{L}}\binom{m-1}{t}\\&\hskip0.5cm -\sum_{t=\mathcal{L}+j-n}^{\mathcal{U}+j-n}\binom{m-1}{t}\\
			&= \sum_{t=max\left\{0,\left\lceil\frac{j-(n-m)}{2}\right\rceil\right\}}^{\left\lceil\frac{m+j}{2}\right\rceil-1}\dbinom{m-1}{t} -\sum_{t=0}^{\left\lfloor\frac{m-j-1}{2}\right\rfloor-1}\dbinom{m-1}{t} \\& \hskip0.5cm- \sum_{t=0}^{\left\lfloor\frac{j-(n-m)}{2}\right\rfloor-1}\dbinom{m-1}{t}.
		\end{aligned}
	\end{equation*}
\end{proof}

To better understand the main theorem, we start by examining a straightforward example. Our goal at this stage is to create a visual representation of weak homomorphisms. Check Figure \ref{Figure-ex1A} for potential weak homomorphisms from $P_4$ to $P_5$ specifically mapping $0$ to $0$. The numbers at the top represent elements of the domain set $V(P_4)$, and those on the left correspond to elements of the image set $V(P_5)$.

The mapping $f_1, f_2 \in \mathrm{WHom}^{0}(P_4,P_5) $ with $f_1(0)=0, f_1(1)=0, f_1(2)=0, f_1(3)=0$ and $f_2(0)=0, f_2(1)=1, f_2(2)=2, f_2(3)=3$ is represented by the dotted line on the top  and black line (see Figure \ref{Figure-ex1C}). 

Figure \ref{Figure-ex1B} illustrates weak homomorphisms using the cubic lattice. Multiple cases need consideration. Initially, when $f(x+1)= f(x)+1$, it corresponds to moving from $(i,j,k)$ to $(i+1,j,k)$. Similarly, when $f(x+1)= f(x)-1$, it corresponds to moving from $(i,j,k)$ to $(i,j+1,k)$. For the remaining cases, where $f(x+1)=f(x)$, the correspondence involves moving from $(i,j,k)$ to $(i,j,k+1)$. Consequently, the mappings $f_1$ and $f_2$ are depicted by the shortest paths from $(0, 0, 0)$ to $(0, 0, 3)$ and $(3, 0, 0)$ in the 0-ladder cubic lattice, respectively.

\begin{figure}[H]
	\begin{minipage}[b]{0.5\textwidth}
		\begin{tikzpicture}[scale=0.65]
			\foreach \x in {6,...,9}
			\draw[fill=black] (\x,5) circle (2pt);
			\foreach \x in {6,...,8}
			\draw [line width=1pt] (\x,5)-- (\x+1,5);
			
			\draw (5.8,5.8) node[anchor=north west] {$0$};
			\draw (6.8,5.8) node[anchor=north west] {$1$};
			\draw (7.8,5.8) node[anchor=north west] {$2$};
			\draw (8.8,5.8) node[anchor=north west] {$3$};

			\foreach \x in {0,...,3}
			\draw [line width=1pt] (5,\x)-- (5,\x+1);

			\draw (4.3,4.3) node[anchor=north west] {$0$};
			\draw (4.3,3.3) node[anchor=north west] {$1$};
			\draw (4.3,2.3) node[anchor=north west] {$2$};
			\draw (4.3,1.3) node[anchor=north west] {$3$};
			\draw (4.3,0.3) node[anchor=north west] {$4$};

			\foreach \x in {6,...,8}
			\draw[dashed] (\x,4) -- (\x+1,4);
			
			\foreach \x in {7,...,8}
			\draw[dashed] (\x,3) -- (\x+1,3);
			
			\draw[dashed] (8,2) -- (9,2);

			\draw[line width=0.3mm] (6,4)-- (7,3);
			\draw[line width=0.3mm] (7,3)-- (8,2);
			\draw[line width=0.3mm] (8,2)-- (9,1);
			\draw[dashed] (7,4)-- (8,3);
			\draw[dashed] (8,3)-- (9,2);
			\draw[line width=0.3mm] (8,4)-- (9,3);
			\draw[line width=0.3mm] (7,3)-- (8,4);
			\draw[dashed] (8,3)-- (9,4);
			\draw[line width=0.3mm] (8,2)-- (9,3);
			
			\foreach \x in {0,...,4}
			\draw[fill=black] (5,\x) circle (2pt);
			
			\foreach \x/\y in {6/4,7/3,8/2,8/4,7/4,8/3,9/4,9/2,9/3,6/4,9/1}
			\draw[fill=black] (\x,\y) circle (2pt);	
			\draw[fill=cyan, line width=0.4mm] (6,4) circle (2.5pt);
			
		\end{tikzpicture}
		\caption{\label{Figure-ex1A}Graphical presentation of domain and image of all possible weak homomorphisms $f : P_4 \to P_5$ where $f(0) =0$.}
	\end{minipage}\hspace{0.5cm}
	\begin{minipage}[b]{0.45\textwidth}
		\begin{tikzpicture}[scale=0.7]	
			\draw[line width=0.3mm,gray] (0,0)-- (3.6,0.6);
			\draw[line width=0.3mm,gray] (-1,0.3)-- (1.4,0.7);
			\draw[line width=0.3mm,gray] (-2,0.6)-- (-0.8,0.8);
			\draw[line width=0.3mm,gray] (1.2,1.2)-- (2.4,1.4);
			
			\draw[line width=0.6mm,black] (0.2,0.5) -- (0.2,1.5);	
			\draw[line width=0.6mm,black] (1.2,0.2) -- (1.2,1.2);
			\draw[line width=0.6mm,black] (2.4,0.4) -- (2.4,1.4);
			
			\draw[dashed] (0,0) -- (-3,0.9);			
			\draw[dashed] (1.2,0.2) -- (-0.8,0.8);			
			\draw[dashed] (2.4,0.4) -- (1.4,0.7);			
			\draw[dashed] (1.2,1.2) -- (0.2,1.5);

			\draw[line width=0.3mm,->,>=stealth,black] (1.5,3) -- (1.4,0.8);
			\draw (0.4,3.8) node[anchor=north west] {$(2,0,1)$};
			
			\draw[fill=cyan, line width=0.4mm] (0,0) circle (2.5pt);
			\draw[fill=gray] (-1,0.3) circle (2pt);
			\draw[fill=gray](-2,0.6) circle (2pt);		
			\draw[fill=black] (-3,0.9) circle (3pt);	
			
			\draw[fill=gray](1.2,0.2) circle (2pt);
			\draw[fill=gray] (0.2,0.5) circle (2pt);
			\draw[fill=black] (-0.8,0.8) circle (3pt);
			
			\draw[fill=gray] (2.4,0.4) circle (2pt);
			\draw[fill=black] (1.4,0.7) circle (3pt);
			
			\draw[fill=black] (3.6,0.6) circle (3pt);
			
			\draw[fill=black] (2.4,1.4) circle (3pt);
			\draw[fill=gray] (1.2,1.2) circle (2pt);
			\draw[fill=black] (0.2,1.5) circle (3pt);
			
			\draw (-2.,1.7) node[anchor=north west] {$(1,0,2)$};
			\draw (-1.1,-0.13) node[anchor=north west] {$(0,0,0)$};
			\draw (2.5,1.5) node[anchor=north west] {$(3,0,0)$};
			\draw (-4.2,1.9) node[anchor=north west] {$(0,0,3)$};
			\draw (1.3,2.3) node[anchor=north west] {$(2,1,0)$};
			\draw (-1.,2.4) node[anchor=north west] {$(1,1,1)$};
		\end{tikzpicture}
		\caption{\label{Figure-ex1B} Cubic lattice presentation of all possible weak homomorphisms $f : P_4 \to P_5$ where $f(0) =0$.}
	\end{minipage}
\end{figure}

\begin{figure}[H]
	\begin{minipage}[b]{0.5\textwidth}
		\begin{tikzpicture}[scale=0.65]
			\foreach \x in {6,...,9}
			\draw[fill=black] (\x,5) circle (2pt);
			\foreach \x in {6,...,8}
			\draw [line width=1pt] (\x,5)-- (\x+1,5);
			
			\draw (5.8,5.8) node[anchor=north west] {$0$};
			\draw (6.8,5.8) node[anchor=north west] {$1$};
			\draw (7.8,5.8) node[anchor=north west] {$2$};
			\draw (8.8,5.8) node[anchor=north west] {$3$};

			\foreach \x in {0,...,3}
			\draw [line width=1pt] (5,\x)-- (5,\x+1);
			
			\draw (4.3,4.3) node[anchor=north west] {$0$};
			\draw (4.3,3.3) node[anchor=north west] {$1$};
			\draw (4.3,2.3) node[anchor=north west] {$2$};
			\draw (4.3,1.3) node[anchor=north west] {$3$};
			\draw (4.3,0.3) node[anchor=north west] {$4$};

			\foreach \x in {6,...,8}
			\draw[dashed] (\x,4) -- (\x+1,4);

			\draw[line width=0.3mm] (6,4)-- (7,3);
			\draw[line width=0.3mm] (7,3)-- (8,2);
			\draw[line width=0.3mm] (8,2)-- (9,1);

			\foreach \x in {0,...,4}
			\draw[fill=black] (5,\x) circle (2pt);
			
			\foreach \x/\y in {6/4,7/3,8/2,8/4,7/4,9/4,6/4,9/1}
			\draw[fill=black] (\x,\y) circle (2pt);	
			\draw[fill=cyan, line width=0.4mm] (6,4) circle (2.5pt);
			
			\draw (9.25,4.5) node[anchor=north west] {$f_1$};
			\draw (7.5,1.5) node[anchor=north west] {$f_2$};
			
		\end{tikzpicture}
		\caption{\label{Figure-ex1C} Graphical presentation of domain and image of $f_1$ and $f_2$.}
	\end{minipage}\hspace{0.5cm}
	\begin{minipage}[b]{0.5\textwidth}
		\begin{tikzpicture}[scale=0.7]	
			\draw[line width=0.3mm,gray] (0,0)-- (3.6,0.6);
			\draw[dashed] (0,0) -- (-3,0.9);			
			
			\draw[fill=cyan, line width=0.4mm] (0,0) circle (2.5pt);
			\draw[fill=gray] (-1,0.3) circle (2pt);
			\draw[fill=gray](-2,0.6) circle (2pt);		
			\draw[fill=black] (-3,0.9) circle (3pt);	
			
			\draw[fill=gray](1.2,0.2) circle (2pt);		
			\draw[fill=gray] (2.4,0.4) circle (2pt);			
			\draw[fill=black] (3.6,0.6) circle (3pt);

			\draw (-1.1,-0.15) node[anchor=north west] {$(0,0,0)$};
			\draw (2.5,1.5) node[anchor=north west] {$(3,0,0)$};
			\draw (-4.2,1.9) node[anchor=north west] {$(0,0,3)$};
			\draw (-2.2,0.5) node[anchor=north west] {$f_1$};
			\draw (1.9,0.4) node[anchor=north west] {$f_2$};
			
		\end{tikzpicture}
		\caption{\label{Figure-ex1D} Cubic lattice presentation of $f_1$ and $f_2$.}
	\end{minipage}
\end{figure}

\noindent The cardinality $|\mathrm{WHom}^{0}(P_4,P_5)|$ is the summation of $M(i, j, k)$ and $M_0(i, j, k)$
where $i + j + k = 3$ (large black points). From Figure \ref{Figure-ex1B}, if $j\leq 0$, we use $M(i, j, k)$, otherwise $M_0(i, j, k)$.

\begin{equation*} \begin{aligned}
	|\mathrm{WHom}^{0}(P_4,P_5)|	&= M(3, 0, 0) + M_0(2, 1, 0) + M(2, 0, 1) + M_0(1, 1, 1) \\
	&\quad + M(1, 0, 2) + M(0, 0, 3)\\
	&= \dbinom{3}{3,0,0} + \left[ \dbinom{3}{2,1,0} - \dbinom{3}{0,3,0}\right]  + \dbinom{3}{2,0,1}\\
	& \quad + \left[ \dbinom{3}{1,1,1} - \dbinom{3}{0,2,1}\right]  
	+\dbinom{3}{1,0,2} + \dbinom{3}{0,0,3}\\
	&= 13.
	\end{aligned}	
	\end{equation*}

\noindent Similar to above example, Figure \ref{Figure-ex2A} visualizes the possible weak homomorphisms of the path $P_4$ to $P_5$ which map $0$ to $1$.  

\begin{figure}[H]
	\begin{minipage}[b]{0.5\textwidth}
		\begin{tikzpicture}[scale=0.65]
			\foreach \x in {6,...,9}
			\draw[fill=black] (\x,5) circle (2pt);
			\foreach \x in {6,...,8}
			\draw [line width=1pt] (\x,5)-- (\x+1,5);
			
			\draw (5.8,5.8) node[anchor=north west] {$0$};
			\draw (6.8,5.8) node[anchor=north west] {$1$};
			\draw (7.8,5.8) node[anchor=north west] {$2$};
			\draw (8.8,5.8) node[anchor=north west] {$3$};

			\foreach \x in {0,...,3}
			\draw [line width=1pt] (5,\x)-- (5,\x+1);

			\draw (4.3,4.3) node[anchor=north west] {$0$};
			\draw (4.3,3.3) node[anchor=north west] {$1$};
			\draw (4.3,2.3) node[anchor=north west] {$2$};
			\draw (4.3,1.3) node[anchor=north west] {$3$};
			\draw (4.3,0.3) node[anchor=north west] {$4$};
			
			\foreach \x in {0,...,4}
			\draw[fill=black] (5,\x) circle (2pt);
			
			\foreach \x in {7,...,8}
			\draw[dashed] (\x,4) -- (\x+1,4);
			
			\foreach \x in {6,...,8}
			\draw[dashed] (\x,3) -- (\x+1,3);
			
			\foreach \x in {7,...,8}
			\draw[dashed] (\x,2) -- (\x+1,2);

			\draw[dashed] (8,1) -- (9,1);
			
			\draw[line width=0.3mm] (6,3)-- (7,2);
			\draw[line width=0.3mm] (6,3)-- (7,4);
			\draw[line width=0.3mm] (7,2)-- (8,3);
			\draw[line width=0.3mm] (7,2)-- (9,0);
			\draw[line width=0.3mm] (8,1)-- (9,2);
			\draw[dashed] (7,3)-- (8,2);
			\draw[dashed] (8,2)-- (9,1);
			\draw[line width=0.3mm] (7,4)-- (8,3);
			\draw[line width=0.3mm] (8,3)-- (9,2);
			\draw[dashed] (8,4)-- (9,3);
			\draw[dashed] (7,3)-- (8,4);
			\draw[line width=0.3mm] (8,3)-- (9,4);
			\draw[dashed] (8,2)-- (9,3);

			\foreach \x/\y in {7/4,7/2,8/1,7/3,8/4,8/3,8/2,6/3,9/4,9/3,9/2,9/1,9/0}
			\draw[fill=black] (\x,\y) circle (2pt);
			
			\draw[fill=cyan,line width=0.4mm] (6,3) circle (2.5pt);
		\end{tikzpicture}
		\caption{\label{Figure-ex2A} Graphical presentation of domain and image of all possible weak homomorphisms $f : P_4 \to P_5$ where $f(0) = 1$}
	\end{minipage}\hspace{0.5cm}
	\begin{minipage}[b]{0.5\textwidth}
		\begin{tikzpicture}[scale=0.8]
			
			\draw[line width=0.3mm,gray] (0,0)-- (3.6,0.6);
			\draw[line width=0.3mm,gray] (-1,0.3)-- (1.4,0.7);
			\draw[line width=0.3mm,gray] (-2,0.6)-- (-0.8,0.8);
			\draw[line width=0.3mm,gray] (0,1)-- (2.4,1.4);
			\draw[line width=0.3mm,gray] (-1,1.3)-- (0.2,1.5);
			
			\draw[line width=0.6mm,black] (0,0) -- (0,1);	
			\draw[line width=0.6mm,black] (-1,0.3) -- (-1,1.3);	
			\draw[line width=0.6mm,black] (-2,0.6) -- (-2,1.6);	
			
			\draw[line width=0.6mm,black] (0.2,0.5) -- (0.2,1.5);	
			\draw[line width=0.6mm,black] (1.2,0.2) -- (1.2,2.2);
			\draw[line width=0.6mm,black] (2.4,0.4) -- (2.4,1.4);
			
			\draw[dashed] (0,0) -- (-3,0.9);
			\draw[dashed] (0,1) -- (-2,1.6);
			\draw[dashed] (1.2,0.2) -- (-0.8,0.8);
			\draw[dashed] (2.4,0.4) -- (1.4,0.7);
			\draw[dashed] (1.2,1.2) -- (0.2,1.5);

			\draw[line width=0.3mm,->,>=stealth,black] (2,3) -- (1.4,0.8);
			\draw (1.1,3.7) node[anchor=north west] {$(2,0,1)$};
			
			\draw[fill=cyan,line width=0.4mm] (0,0) circle (2.5pt);
			\draw[fill=gray] (-1,0.3) circle (2pt);
			\draw[fill=gray](-2,0.6) circle (2pt);		
			\draw[fill=black] (-3,0.9) circle (3pt);	
			
			\draw[fill=gray] (0,1) circle (2pt);
			\draw[fill=gray] (-1,1.3) circle (2pt);
			\draw[fill=black] (-2,1.6) circle (3pt);	
			
			\draw[fill=gray](1.2,0.2) circle (2pt);
			\draw[fill=gray] (0.2,0.5) circle (2pt);
			\draw[fill=black] (-0.8,0.8) circle (3pt);
			
			\draw[fill=gray] (2.4,0.4) circle (2pt);
			\draw[fill=black] (1.4,0.7) circle (3pt);
			
			\draw[fill=gray] (1.2,1.2) circle (2pt);
			\draw[fill=black] (0.2,1.5) circle (3pt);
			
			\draw[fill=black] (2.4,1.4) circle (3pt);
			
			\draw[fill=black] (1.2,2.2) circle (3pt);
			\draw[fill=black](3.6,0.6) circle (3pt);		
			
			\draw (-1.,-0.15) node[anchor=north west] {$(0,0,0)$};
			\draw (-4.0,1.8) node[anchor=north west] {$(0,0,3)$};
			\draw (1.6,2.2) node[anchor=north west] {$(2,1,0)$};
			\draw (-0.8,2.3) node[anchor=north west] {$(1,1,1)$};
			\draw (-3.0,2.5) node[anchor=north west] {$(0,1,2)$};	
			\draw (0.,3.) node[anchor=north west] {$(1,2,0)$};	
			\draw (-2.0,3.8) node[anchor=north west] {$(1,0,2)$};	
			\draw (2.6,0.5) node[anchor=north west] {$(3,0,0)$};
			\draw[line width=0.3mm,->,>=stealth,black] (-1,3) -- (-0.8,0.9);
			
		\end{tikzpicture}
		\caption{\label{Figure-ex2B} Cubic lattice presentation of all possible weak homomorphisms $f : P_4 \to P_5$ where $f(0) = 1$}
	\end{minipage}
\end{figure}

\noindent The cardinality $|\mathrm{WHom}^{1}(P_4,P_5)|$ is the summation of $M(i, j, k)$ and $M_1(i, j, k)$
where $i + j + k = 3$ (large black points). From Figure \ref{Figure-ex2B}, if $j\leq 1$ , we use $M(i, j, k)$, otherwise $M_1(i, j, k)$.

\begin{equation*} \begin{aligned}
	|\mathrm{WHom}^{1}(P_4,&P_5)|	  = M(2, 1, 0) + M_1(1, 2, 0) + M(2, 0, 1) + M(1, 1, 1)  \\
	& \quad + M(1, 0, 2) + M(0, 1, 2) + M(0, 0, 3)+ M(3, 0, 0)\\
	&= \dbinom{3}{2,1,0} + \left[ \dbinom{3}{1,2,0} - \dbinom{3}{0,3,0}\right]  + \dbinom{3}{2,0,1} + \dbinom{3}{1,1,1}  \\
	& \quad+ \dbinom{3}{1,0,2} + \dbinom{3}{0,1,2} + \dbinom{3}{0,0,3}+ \dbinom{3}{3,0,0}\\
	&= 22.
	\end{aligned}	
	\end{equation*}

In this paper, our interest lies in determining the count of weak homomorphisms from paths to rectangular grid graphs, denoted as |$\mathrm{WHom}(P_m,P_n\square P_k)$|, which provides a solution to the problem concerning the number of partial walks of $m$ vertices within the rectangular grid graphs $P_n\square P_k$.

\section{The Number of Weak Homomorphisms from Paths to Paths that map $0$ to $j$}
\label{2}
In this section, we present the formula for determining the count of weak homomorphisms from paths $P_m$ to $P_n$, where $0$ is mapped to $j$. We represent the set of weak homomorphisms from $P_m$ to $P_n$, with the mapping of $0$ to $j$, as WHom$^j$($P_m, P_n$).

\begin{theorem}\label{mainthm1}
	Let $m,n$ be positive integers and $j$ a non-negative integer such that $m\leq n$ and $j < n$.  Then
	
	\begin{equation*}
	\begin{aligned}
		|\mathrm{WHom}^{j}&(P_m,P_n)| =\sum_{t=j+1}^{j+\left\lfloor \frac{m-j-1}{2}\right\rfloor}\sum_{s=t-j}^{m-1-t}\left[\tbinom{m-1}{s,t,m-1-s-t} - \tbinom{m-1}{t-j-1,s+j+1,m-1-s-t}\right]\\
		&+\sum_{t=max\{j-n+m+1,0\}}^{j}\sum_{s=0}^{m-1-t}\tbinom{m-1}{s,t,m-1-s-t}+\sum_{t=0}^{j-n+m}\sum_{s=0}^{n-j-1}\tbinom{m-1}{s,t,m-1-s-t}\\
		&+\sum_{t=n-j-1+1}^{n-j-1+\left\lfloor\frac{j-n+m}{2}\right\rfloor}\sum_{s=t-(n-j-1)}^{m-1-t}\left[\tbinom{m-1}{s,t,m-1-s-t} - \tbinom{m-1}{t-n+j,s+n-j,m-1-s-t}\right].
	\end{aligned}	
\end{equation*}
	
\end{theorem}

\begin{proof}
	To find $|\mathrm{WHom}^{j}(P_m, P_n)|$, we count the number of shortest paths from the point $(0, 0, 0)$ to any point $(i_0, j_0, k_0)$, where $i_0+j_0+k_0=m-1$ in the $j$-ladder cubic lattice. Consider the following three different cases correspond to the value of $j_0$: \\
	
	\textbf{Case 1: } $j_0 > j$. For each $j_0=j+t$, there are $\sum_{i_0=t}^{m-j-1-t}M_j(i_0,j_0,k_0)$ shortest paths.\\
	
	\noindent Since $t \leq \frac{m-j-1}{2}$, we obtain \\	
	\begin{equation*}
	\begin{aligned}
		\sum_{t=1}^{\left\lfloor \frac{m-j-1}{2}\right\rfloor}\sum_{i_0=t}^{m-j-1-t}&M_j(i_0,j+t,k_0)
		=\sum_{t=j+1}^{j+\left\lfloor \frac{m-j-1}{2}\right\rfloor}\sum_{s=t-j}^{m-1-t}M_j(s,t,m-1-s-t)\\
		&=\sum_{t=j+1}^{j+\left\lfloor \frac{m-j-1}{2}\right\rfloor}\sum_{s=t-j}^{m-1-t}\left[\tbinom{m-1}{s,t,m-1-s-t} - \tbinom{m-1}{t-j-1,s+j+1,m-1-s-t}\right].
	\end{aligned}	
\end{equation*}
	
	\begin{figure}[H]\begin{center}
			\begin{tikzpicture}[scale=0.3]
				\def\m{13}
				\def\n{15}
				\def\r{4}
				\pgfmathsetmacro{\mnr}{\m-\n+\r}
				\pgfmathsetmacro{\nr}{\n-\r-1}
				\pgfmathsetmacro{\mr}{\m-\r-1}
				\pgfmathsetmacro{\hmr}{\mr/2}
				\pgfmathsetmacro{\hnr}{\nr/2}
				\pgfmathsetmacro{\hmnr}{\mnr/2}
				
				\draw (-1.5,-5) node[anchor=north] {$(0,0)$};
				\draw (-1.5,\r) node[anchor=south] {$(0,j)$};		
				\draw (\mr,-5) node[anchor=north] {$(m-j-1,0)$};
				\draw (16.5,-4) node[anchor=north] {$(n-j-1,0)$};
				\draw (-4.5,1.2) node[anchor=north] {$(0,j-n+m)$};

				\draw[fill=black] (0,-5) circle(1mm);
				\draw[fill=black] (\mr,-5) circle(1mm);
				\draw[fill=black] (0,\r) circle(1mm);
				\draw[fill=black] (12,-5) circle(1mm);
				\draw[fill=black] (0,0) circle(1mm);

				\foreach \x in {0,1,...,{\hmr}} 
				{\draw (\x, 0) -- (\x,{\r+\x});
					\draw ({\x+\hmr}, 0) -- ({\x+\hmr},{\r+\hmr-\x});}
				
				\foreach \x in {1,2,3,4} 
				{\draw ({\x+8}, -5) -- ({\x+8},{\r-\x});}
				
				\foreach \x in {0,...,8} 
				{\draw (\x, -5) -- ({\x},0);}
				
				\foreach \x in {1,2} 
				{\draw (\x+12, -5+\x) -- ({\x+12},-\x);}
				
				\foreach \y in {0,...,\r}{			
					\draw (0, \y) -- (\mr,\y);}
				
				\foreach \y in {1,...,\hmr}			
				\draw (\y, \r+\y) -- (\mr-\y, \r+\y);
				
				\foreach \y in {1,2,3}			
				\draw (0, \y-6) -- (11+\y, \y-6);
				
				\foreach \y in {4,...,9}			
				\draw (0, \y-6) -- (18-\y, \y-6);
				
				\draw[line width=0.5mm,->,>=stealth,blue] (5,7)--(3,7);
				\draw[line width=0.5mm,->,>=stealth,blue] (6,6)--(2,6);
				\draw[line width=0.5mm,->,>=stealth,blue] (7,5)--(1,5);
				
				\draw[line width=0.4mm, black] (0,-5) -- (12,-5);
				\draw[line width=0.4mm, black] (0,-5) -- (0,\r);
				\draw[line width=0.4mm, black] (0,\r) -- (\mr,\r);
				\draw[line width=0.4mm, black] (\mr,\r) -- (\mr,-5);
				\draw[line width=0.4mm, black] (12,0) -- (0,0);
				\draw[line width=0.4mm, black] (12,0) -- (12,-5);

				\foreach \y in {1,...,\hmr}
				\draw[fill=magenta] (\mr-\y, \r+\y) circle(1mm) ;
				
		\end{tikzpicture}\end{center}
		\caption{\label{FigureP1} Points $(i_0,j_0)$ where $j_0>j$. }
	\end{figure}
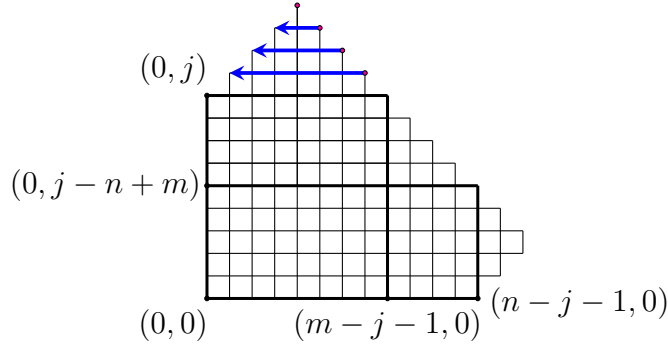

	\textbf{Case 2: } $j-n+m< j_0\leq j$ and $i_0 < n-j-1$. For each $j_0=j-t$, there are $\sum_{i_0=0}^{m-j-1+t}M(i_0,j_0,k_0)$ shortest paths.\\
	
	\noindent Since $t < n-m$, we obtain 
	
	\begin{equation*}
	\begin{aligned}
		\sum_{t=0}^{n-m-1}\sum_{i_0=0}^{m-j-1+t}M(i_0,j-t,k_0)&=\sum_{t=max\{j-n+m+1,0\}}^{j}\sum_{s=0}^{m-1-t}M(s,t,m-1-s-t)\\
		&=\sum_{t=max\{j-n+m+1,0\}}^{j}\sum_{s=0}^{m-1-t}\tbinom{m-1}{s,t,m-1-s-t}.
	\end{aligned}	
\end{equation*}
	
	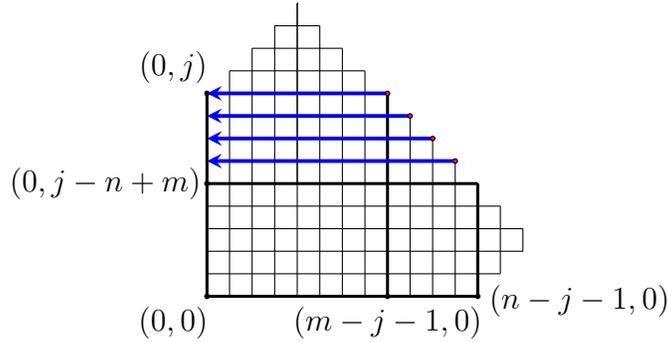
\begin{figure}[H]\begin{center}
			\begin{tikzpicture}[scale=0.3]
				\def\m{13}
				\def\n{15}
				\def\r{4}
				\pgfmathsetmacro{\mnr}{\m-\n+\r}
				\pgfmathsetmacro{\nr}{\n-\r-1}
				\pgfmathsetmacro{\mr}{\m-\r-1}
				\pgfmathsetmacro{\hmr}{\mr/2}
				\pgfmathsetmacro{\hnr}{\nr/2}
				\pgfmathsetmacro{\hmnr}{\mnr/2}
				

				\draw (-1.5,-5) node[anchor=north] {$(0,0)$};
				\draw (-1.5,\r) node[anchor=south] {$(0,j)$};		
				\draw (\mr,-5) node[anchor=north] {$(m-j-1,0)$};
				\draw (16.5,-4) node[anchor=north] {$(n-j-1,0)$};
				\draw (-4.5,1.2) node[anchor=north] {$(0,j-n+m)$};
				
				\draw[fill=black] (0,-5) circle(1mm);
				\draw[fill=black] (\mr,-5) circle(1mm);
				\draw[fill=black] (0,\r) circle(1mm);
				\draw[fill=black] (12,-5) circle(1mm);
				\draw[fill=black] (0,0) circle(1mm);

				\foreach \x in {0,1,...,{\hmr}} 
				{\draw (\x, 0) -- (\x,{\r+\x});
					\draw ({\x+\hmr}, 0) -- ({\x+\hmr},{\r+\hmr-\x});}
				
				\foreach \x in {1,2,3,4} 
				{\draw ({\x+8}, -5) -- ({\x+8},{\r-\x});}
				
				\foreach \x in {0,...,8} 
				{\draw (\x, -5) -- ({\x},0);}
				
				\foreach \x in {1,2} 
				{\draw (\x+12, -5+\x) -- ({\x+12},-\x);}
				
				\foreach \y in {0,...,\r}{			
					\draw (0, \y) -- (\mr,\y);}
				
				\foreach \y in {1,...,\hmr}			
				\draw (\y, \r+\y) -- (\mr-\y, \r+\y);
				
				\foreach \y in {1,2,3}			
				\draw (0, \y-6) -- (11+\y, \y-6);
				
				\foreach \y in {4,...,9}			
				\draw (0, \y-6) -- (18-\y, \y-6);

				\draw[line width=0.4mm, black] (0,-5) -- (12,-5);
				\draw[line width=0.4mm, black] (0,-5) -- (0,\r);
				\draw[line width=0.4mm, black] (0,\r) -- (\mr,\r);
				\draw[line width=0.4mm, black] (\mr,\r) -- (\mr,-5);
				\draw[line width=0.4mm, black] (12,0) -- (0,0);
				\draw[line width=0.4mm, black] (12,0) -- (12,-5);
				
				\draw[line width=0.5mm,->,>=stealth,blue] (8,4)--(0,4);
				\draw[line width=0.5mm,->,>=stealth,blue] (9,3)--(0,3);
				\draw[line width=0.5mm,->,>=stealth,blue] (10,2)--(0,2);
				\draw[line width=0.5mm,->,>=stealth,blue] (11,1)--(0,1);
				
				\foreach \y in {1,...,4}
				\draw[fill=magenta] (7+\y, \r-\y+1) circle(1mm) ;
				
		\end{tikzpicture}\end{center}
		\caption{\label{FigureP2} Points $(i_0,j_0)$ where $j-n+m<j_0\leq j$.}
	\end{figure}

	\textbf{Case 3: } $j_0\leq j-n+m\leq j$ and $i_0 \leq n-j-1$. For each $i_0, j_0$, there are $M(i_0,j_0,k_0)$ shortest paths. \\
	
	\noindent We obtain  
	
	\begin{equation*}
		\begin{aligned}
			\sum_{j_0=0}^{j-n+m}\sum_{i_0=0}^{n-j-1}M(i_0,j_0,k_0)&=\sum_{t=0}^{j-n+m}\sum_{s=0}^{n-j-1}M(s,t,m-1-s-t)\\
			&=\sum_{t=0}^{j-n+m}\sum_{s=0}^{n-j-1}\tbinom{m-1}{s,t,m-1-s-t}.
		\end{aligned}	
	\end{equation*}
	
	\begin{figure}[H]\begin{center}
			\begin{tikzpicture}[scale=0.3]
				\def\m{13}
				\def\n{15}
				\def\r{4}
				\pgfmathsetmacro{\mnr}{\m-\n+\r}
				\pgfmathsetmacro{\nr}{\n-\r-1}
				\pgfmathsetmacro{\mr}{\m-\r-1}
				\pgfmathsetmacro{\hmr}{\mr/2}
				\pgfmathsetmacro{\hnr}{\nr/2}
				\pgfmathsetmacro{\hmnr}{\mnr/2}
				

				\draw (-1.5,-5) node[anchor=north] {$(0,0)$};
				\draw (-1.5,\r) node[anchor=south] {$(0,j)$};		
				\draw (\mr,-5) node[anchor=north] {$(m-j-1,0)$};
				\draw (16.5,-4) node[anchor=north] {$(n-j-1,0)$};
				\draw (-4.5,1.2) node[anchor=north] {$(0,j-n+m)$};

				\draw[fill=black] (0,-5) circle(1mm);
				\draw[fill=black] (\mr,-5) circle(1mm);
				\draw[fill=black] (0,\r) circle(1mm);
				\draw[fill=black] (12,-5) circle(1mm);
				\draw[fill=black] (0,0) circle(1mm);

				\foreach \x in {0,1,...,{\hmr}} 
				{\draw (\x, 0) -- (\x,{\r+\x});
					\draw ({\x+\hmr}, 0) -- ({\x+\hmr},{\r+\hmr-\x});}
				
				\foreach \x in {1,2,3,4} 
				{\draw ({\x+8}, -5) -- ({\x+8},{\r-\x});}
				
				\foreach \x in {0,...,8} 
				{\draw (\x, -5) -- ({\x},0);}
				
				\foreach \x in {1,2} 
				{\draw (\x+12, -5+\x) -- ({\x+12},-\x);}
				
				\foreach \y in {0,...,\r}{			
					\draw (0, \y) -- (\mr,\y);}
				
				\foreach \y in {1,...,\hmr}			
				\draw (\y, \r+\y) -- (\mr-\y, \r+\y);
				
				\foreach \y in {1,2,3}			
				\draw (0, \y-6) -- (11+\y, \y-6);
				
				\foreach \y in {4,...,9}			
				\draw (0, \y-6) -- (18-\y, \y-6);

				\draw[line width=0.4mm, black] (0,-5) -- (12,-5);
				\draw[line width=0.4mm, black] (0,-5) -- (0,\r);
				\draw[line width=0.4mm, black] (0,\r) -- (\mr,\r);
				\draw[line width=0.4mm, black] (\mr,\r) -- (\mr,-5);
				\draw[line width=0.4mm, black] (12,0) -- (0,0);
				\draw[line width=0.4mm, black] (12,0) -- (12,-5);

				\draw[line width=0.5mm,->,>=stealth,blue] (12,-1)--(0,-1);
				\draw[line width=0.5mm,->,>=stealth,blue] (12,-2)--(0,-2);
				\draw[line width=0.5mm,->,>=stealth,blue] (12,-3)--(0,-3);
				\draw[line width=0.5mm,->,>=stealth,blue] (12,-4)--(0,-4);
				\draw[line width=0.5mm,->,>=stealth,blue] (12,-5)--(0,-5);
				\draw[line width=0.5mm,->,>=stealth,blue] (12,0)--(0,0);
				
				\foreach \y in {1,...,6}
				\draw[fill=magenta] (12, -6+\y) circle(1mm) ;
				
		\end{tikzpicture}\end{center}
		\caption{\label{FigureP3} Points $(i_0,j_0)$ where $j_0<j-n+m$ and $i_0\leq n-j-1$.}
	\end{figure}
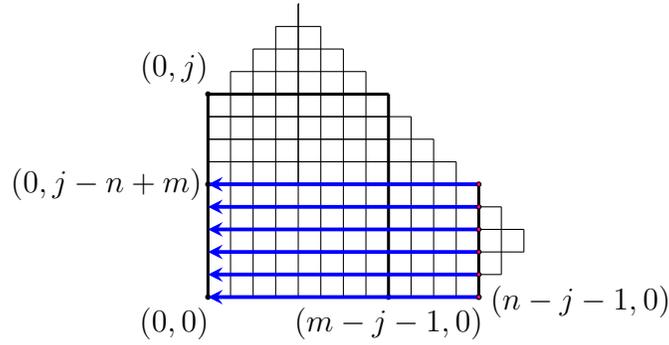

	\textbf{Case 4: } $j_0\leq j$ and $i_0 > n-j-1$. For each $i_0=n-j-1+t$, there are $\sum_{j_0=t}^{j-n+m-t}M_{n-j-1}(j_0,i_0,k_0)$ shortest paths. This can be obtained by flipping the cubic lattice diagonally. 
	
	\noindent	Since $t \leq \frac{j-n+m}{2}$, we obtain
	\begin{equation*}
		\begin{aligned}
			\sum_{t=1}^{\left\lfloor\frac{j-n+m}{2}\right\rfloor}&\sum_{j_0=t}^{j-n+m-t}M_{n-j-1}(j_0,n-j-1+t,k_0)\\&=\sum_{t=n-j-1+1}^{n-j-1+\left\lfloor\frac{j-n+m}{2}\right\rfloor}\sum_{s=t-(n-j-1)}^{m-1-t}M_{n-j-1}(s,t,m-1-s-t)\\
			&=\sum_{t=n-j-1+1}^{n-j-1+\left\lfloor\frac{j-n+m}{2}\right\rfloor}\sum_{s=t-(n-j-1)}^{m-1-t}\left[\tbinom{m-1}{s,t,m-1-s-t} - \tbinom{m-1}{t-n+j,s+n-j,m-1-s-t}\right].
		\end{aligned}	
	\end{equation*}

	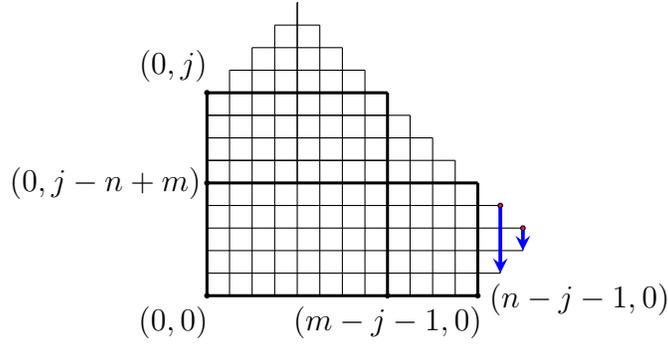
\begin{figure}[H]\begin{center}
			\begin{tikzpicture}[scale=0.3]
				\def\m{13}
				\def\n{15}
				\def\r{4}
				\pgfmathsetmacro{\mnr}{\m-\n+\r}
				\pgfmathsetmacro{\nr}{\n-\r-1}
				\pgfmathsetmacro{\mr}{\m-\r-1}
				\pgfmathsetmacro{\hmr}{\mr/2}
				\pgfmathsetmacro{\hnr}{\nr/2}
				\pgfmathsetmacro{\hmnr}{\mnr/2}
				
				
				\draw (-1.5,-5) node[anchor=north] {$(0,0)$};
				\draw (-1.5,\r) node[anchor=south] {$(0,j)$};		
				\draw (\mr,-5) node[anchor=north] {$(m-j-1,0)$};
				\draw (16.5,-4) node[anchor=north] {$(n-j-1,0)$};
				\draw (-4.5,1.2) node[anchor=north] {$(0,j-n+m)$};
				
				\draw[fill=black] (0,-5) circle(1mm);
				\draw[fill=black] (\mr,-5) circle(1mm);
				\draw[fill=black] (0,\r) circle(1mm);
				\draw[fill=black] (12,-5) circle(1mm);
				\draw[fill=black] (0,0) circle(1mm);

				\foreach \x in {0,1,...,{\hmr}} 
				{\draw (\x, 0) -- (\x,{\r+\x});
					\draw ({\x+\hmr}, 0) -- ({\x+\hmr},{\r+\hmr-\x});}
				
				\foreach \x in {1,2,3,4} 
				{\draw ({\x+8}, -5) -- ({\x+8},{\r-\x});}
				
				\foreach \x in {0,...,8} 
				{\draw (\x, -5) -- ({\x},0);}
				
				\foreach \x in {1,2} 
				{\draw (\x+12, -5+\x) -- ({\x+12},-\x);}
				
				\foreach \y in {0,...,\r}{			
					\draw (0, \y) -- (\mr,\y);}
				
				\foreach \y in {1,...,\hmr}			
				\draw (\y, \r+\y) -- (\mr-\y, \r+\y);
				
				\foreach \y in {1,2,3}			
				\draw (0, \y-6) -- (11+\y, \y-6);
				
				\foreach \y in {4,...,9}			
				\draw (0, \y-6) -- (18-\y, \y-6);

				\draw[line width=0.4mm, black] (0,-5) -- (12,-5);
				\draw[line width=0.4mm, black] (0,-5) -- (0,\r);
				\draw[line width=0.4mm, black] (0,\r) -- (\mr,\r);
				\draw[line width=0.4mm, black] (\mr,\r) -- (\mr,-5);
				\draw[line width=0.4mm, black] (12,0) -- (0,0);
				\draw[line width=0.4mm, black] (12,0) -- (12,-5);

				\draw[line width=0.5mm,->,>=stealth,blue] (13,-1)--(13,-4);
				\draw[line width=0.5mm,->,>=stealth,blue] (14,-2)--(14,-3);
				
				\draw[fill=magenta] (13,-1) circle(1mm) ;
				\draw[fill=magenta] (14,-2) circle(1mm) ;
				
		\end{tikzpicture}\end{center}
		\caption{\label{FigureP4} Points $(i_0,j_0)$ where $j_0<j-n+m$ and $i_0> n-j-1$.}
	\end{figure}
	
	Adding up over all cases, $|\mathrm{WHom}^{j}(P_m,P_n)|$ is as desired. 
	
\end{proof}

For convenience, we compute $|\mathrm{Hom}^j(P_m,P_n)|$ and $|\mathrm{WHom}^j(P_m,P_n)|$ for $2\leq m\leq n\leq8$. The results are presented in Tables \ref{tab3} and \ref{tab4}, respectively.

\begin{table}[H]
	\caption{\label{tab3} Numbers of weak homomorphisms $f:P_m \rightarrow P_n$ where $f(0)=j$ for $2\leq m\leq n\leq 9.$}
	\begin{minipage}{0.5\textwidth}
		\begin{tabular}{c|c|c|c|c|c|c|c|c|}
			\multicolumn{2}{c}{ } & \multicolumn{7}{|c|}{$n$}\\\hline
			$m$&$j$&2&3&4&5&6&7&8\\\hline
			\multirow{4}{*}{2}
			&0&2&2&2&2&2&2&2\\\cline{2-9}
			&1&&3&3&3&3&3&3\\\cline{2-9}
			&2&&&&3&3&3&3\\\cline{2-9}
			&3&&&&&&3&3\\\hline
			\multirow{4}{*}{3}
			&0&&5&5&5&5&5&5\\\cline{2-9}
			&1&&7&8&8&8&8&8\\\cline{2-9}
			&2&& & &9&9&9&9\\\cline{2-9}
			&3&& & & & &9&9\\\hline
			\multirow{4}{*}{4}
			&0&&&13&13&13&13&13\\\cline{2-9}
			&1&&&21&22&22&22&22\\\cline{2-9}
			&2&&&  &25&26&26&26\\\cline{2-9}
			&3&&&  &  &  &27&27\\\hline
			\multirow{4}{*}{5}
			&0&&&&35&35&35&35\\\cline{2-9}
			&1&&&&60&61&61&61\\\cline{2-9}
			&2&&&&69&74&75&75\\\cline{2-9}
			&3&&&&  &  &79&80\\\hline
			
			\multirow{4}{*}{6}
			&0&&&&&96&96&96\\\cline{2-9}
			&1&&&&&170&171&171\\\cline{2-9}
			&2&&&&&209&215&216\\\cline{2-9}
			&3&&&&&   &229&235\\\hline
			\multirow{4}{*}{7}
			&0&&&&&&267&267\\\cline{2-9}
			&1&&&&&&482&483\\\cline{2-9}
			&2&&&&&&615&622\\\cline{2-9}
			&3&&&&&&659&686\\\hline
			\multirow{4}{*}{8}
			&0&&&&&&&750\\\cline{2-9}
			&1&&&&&&&1372\\\cline{2-9}
			&2&&&&&&&1791\\\cline{2-9}
			&3&&&&&&&1994\\\hline
		\end{tabular}
	\end{minipage}
\end{table}

\begin{table}[H]
	\caption{\label{tab4} Numbers of homomorphisms $f:P_m \rightarrow P_n$ where $f(0)=j$ for $2\leq m\leq n\leq 9.$}
	\begin{minipage}{0.5\textwidth}
		\begin{tabular}{c|c|c|c|c|c|c|c|c|}
			\multicolumn{2}{c}{ } & \multicolumn{7}{|c|}{$n$}\\\hline
			$m$&$j$&2&3&4&5&6&7&8\\\hline
			\multirow{4}{*}{2}&0&1&1&1&1&1&1&1\\\cline{2-9}
			&1&&2&2&2&2&2&2\\\cline{2-9}
			&2&&&&2&2&2&2\\\cline{2-9}
			&3&&&&&&2&2\\\hline
			\multirow{4}{*}{3}&0&&2&2&2&2&2&2\\\cline{2-9}
			&1&&2&3&3&3&3&3\\\cline{2-9}
			&2&&&&4&4&4&4\\\cline{2-9}
			&3&&&&&&4&4\\\hline
			\multirow{4}{*}{4}&0&&&3&3&3&3&3\\\cline{2-9}
			&1&&&5&6&6&6&6\\\cline{2-9}
			&2&&&&6&7&7&7\\\cline{2-9}
			&3&&&&&&8&8\\\hline
			\multirow{4}{*}{5}&0&&&&6&6&6&6\\\cline{2-9}
			&1&&&&9&10&10&10\\\cline{2-9}
			&2&&&&12&13&14&14\\\cline{2-9}
			&3&&&&&&14&15\\\hline
			
			\multirow{4}{*}{6}&0&&&&&10&10&10\\\cline{2-9}
			&1&&&&&19&20&20\\\cline{2-9}
			&2&&&&&23&24&25\\\cline{2-9}
			&3&&&&&&28&29\\\hline
			\multirow{4}{*}{7}&0&&&&&&20&20\\\cline{2-9}
			&1&&&&&&34&35\\\cline{2-9}
			&2&&&&&&48&49\\\cline{2-9}
			&3&&&&&&48&54\\\hline
			\multirow{4}{*}{8}&0&&&&&&&35\\\cline{2-9}
			&1&&&&&&&69\\\cline{2-9}
			&2&&&&&&&89\\\cline{2-9}
			&3&&&&&&&103\\\hline
		\end{tabular}
	\end{minipage}
	
\end{table}

 \section{The Number of Weak Homomorphisms from Paths to Grid Graphs}
 \label{mainsec}

In this section, we present the formulas for determining the count of weak homomorphisms from paths $P_m$ to rectangular grid graphs $P_n \square P_k$. We represent the set of weak homomorphisms from $P_m$ to $P_n \square P_k$, mapping $0$ to $(i, j)$, as WHom$^{ij}(P_m, P_n \square P_k)$. From the symmetry of $P_n \square P_k$, we deduce the following lemma:

\begin{lemma}\label{Lemma4}Let $i$ and $n$ be integers such that $0\leq j <n$, and let $m>2$ be a positive integer.
	\begin{enumerate}
		\item $|\mathrm{WHom}^{ij}(P_m,P_n\square P_k)|=| \mathrm{WHom}^{(n-i-1)j}(P_m,P_n\square P_k)| $\\\hspace*{4.17cm}$=| \mathrm{WHom}^{i(k-j-1)}(P_m,P_n\square P_k)|$\\\hspace*{4.17cm}$=| \mathrm{WHom}^{(n-i-1)(k-j-1)}(P_m,P_n\square P_k)|$,\\\hspace*{1.5cm} for all $i\in\{0,1,\dots,n-1\}$ and $j\in\{0,1,\dots,k-1\}$. 	
		\item $| \mathrm{WHom}(P_m,P_{2n}\square P_{2k})|=4\sum_{i=0}^{n-1}\sum_{j=0}^{k-1}| \mathrm{WHom}^{ij}(P_m,P_{2n}\square P_{2k})| $. 	
		\item $| \mathrm{WHom}(P_m,P_{2n+1}\square P_{2k})|=4\sum_{i=0}^{n-1}\sum_{j=0}^{k-1}| \mathrm{WHom}^{ij}(P_m,P_{2n+1}\square P_{2k})|$ \\\hspace*{5cm}
		$+2\sum_{j=0}^{k-1}| \mathrm{WHom}^{nj}(P_m,P_{2n+1}\square P_{2k})|$. 
		\item $| \mathrm{WHom}(P_m,P_{2n}\square P_{2k+1})|=4\sum_{i=0}^{n-1}\sum_{j=0}^{k-1}| \mathrm{WHom}^{ij}(P_m,P_{2n}\square P_{2k+1})|$ \\\hspace*{5cm} $+2\sum_{i=0}^{n-1}| \mathrm{WHom}^{ik}(P_m,P_{2n}\square P_{2k+1})|$.  	
		\item $| \mathrm{WHom}(P_m,P_{2n+1}\square P_{2k+1})|$\\\hspace*{4.cm}$=4\sum_{i=0}^{n-1}\sum_{j=0}^{k-1}| \mathrm{WHom}^{ij}(P_m,P_{2n+1}\square P_{2k+1})| $\\\hspace*{4.5cm}$+2\sum_{j=0}^{k-1}| \mathrm{WHom}^{nj}(P_m,P_{2n+1}\square P_{2k+1})|$\\\hspace*{4.5cm}$+2\sum_{i=0}^{n-1}| \mathrm{WHom}^{ik}(P_m,P_{2n+1}\square P_{2k+1})|$\\\hspace*{4.5cm}$+|\mathrm{WHom}^{nk}(P_m,P_{2n+1}\square P_{2k+1})|$.
	\end{enumerate}	
\end{lemma}

\begin{example}
	$\mathrm{WHom}^{00}(P_4,P_4 \square P_5) = 43$.
\end{example}
	
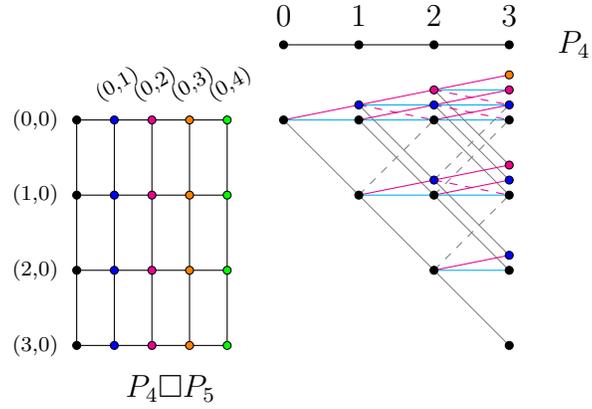
\begin{figure}[H]\begin{center}
			\begin{tikzpicture}[scale=0.5]
				\draw (3.5*2,4+4) node[anchor=west] {$P_4$};
				\draw (-3,-0.5) node[anchor=north] {$P_4\square P_5$};
				
				\draw (-1.7-0.5,6.5) node[anchor=west,rotate=30] {\scriptsize (0,4)};
				\draw (-2.7-0.5,6.5) node[anchor=west,rotate=30] {\scriptsize (0,3)};
				\draw (-3.7-0.5,6.5) node[anchor=west,rotate=30] {\scriptsize (0,2)};
				\draw (-4.7-0.5,6.5) node[anchor=west,rotate=30] {\scriptsize (0,1)};
				
				\draw (-7-0.5,3*2) node[anchor=west] {\scriptsize (0,0)};
				\draw (-7-0.5,2*2) node[anchor=west] {\scriptsize (1,0)};
				\draw (-7-0.5,1*2) node[anchor=west] {\scriptsize (2,0)};
				\draw (-7-0.5,0*2) node[anchor=west] {\scriptsize (3,0)};
				
				\foreach \x in {0,...,3}
				{\draw[fill=black] (\x*2,4+4) circle (3pt);
					\draw (\x*2,4.2+4) node[anchor=south] {\x};}
				\foreach \x in {0,...,2}
				{\draw (\x*2,4+4)-- (2*\x+2,4+4);}
				
				\foreach \y in {0,...,3}
				\draw (-1-0.5,2*\y)-- (-5-0.5,2*\y);
				\foreach \x in {-1,...,-5}
				\draw (\x-0.5,0)-- (\x-0.5,3+3);				
				
				\foreach \y in {0,...,3}
				{\draw[fill=green] (-1-0.5,{-1*2*\y+3+3}) circle (3pt);
					\draw[fill=orange] (-2-0.5,{-1*2*\y+3+3}) circle (3pt);
					\draw[fill=magenta] (-3-0.5,{-1*2*\y+3+3}) circle (3pt);
					\draw[fill=blue] (-4-0.5,{-1*2*\y+3+3}) circle (3pt);
					\draw[fill=black] (-5-0.5,{-1*2*\y+3+3}) circle (3pt);}
				
				\foreach \x/\y/\z/\u in {0/3/3/0,1/3/3/1,2/3/3/2}
				\draw[gray] (2*\x,2*\y) -- (2*\z,2*\u);			
				\foreach \x/\y/\z/\u in {1/3/3/1,2/3/3/2}
				\draw[gray] (2*\x,2*\y+0.4) -- (2*\z,2*\u+0.4);			
				\foreach \x/\y/\z/\u in {2/3/3/2}
				\draw[gray] (2*\x,2*\y+0.8) -- (2*\z,2*\u+0.8);	
				
				\foreach \x/\y/\z/\u in {1/2/2/3,2/2/3/3,2/1/3/2}
				\draw[gray,dashed] (2*\x,2*\y) -- (2*\z,2*\u);			
				\foreach \x/\y/\z/\u in {2/2/3/3}
				\draw[gray,dashed] (2*\x,2*\y+0.4) -- (2*\z,2*\u+0.4);			
				
				\foreach \x/\y/\z/\u in {0/3/1/3,1/3/2/3,2/3/3/3,1/2/2/2,2/2/3/2,2/1/3/1}
				\draw[magenta] (2*\x,2*\y) -- (2*\z,2*\u+0.4);								
				\foreach \x/\y/\z/\u in {1/3/2/3,2/3/3/3,2/2/3/2}
				\draw[magenta] (2*\x,2*\y+0.4) -- (2*\z,2*\u+0.8);									
				\foreach \x/\y/\z/\u in {2/3/3/3}
				\draw[magenta] (2*\x,2*\y+0.8) -- (2*\z,2*\u+1.2);	
				
				\foreach \x/\y/\z/\u in {1/3/2/3,2/3/3/3,2/2/3/2}
				\draw[magenta,dashed] (2*\x,2*\y+0.4) -- (2*\z,2*\u);						
				\foreach \x/\y/\z/\u in {2/3/3/3}
				\draw[magenta,dashed] (2*\x,2*\y+0.8) -- (2*\z,2*\u+0.4);	
				
				\foreach \x/\y/\z/\u in {0/3/3/3,1/2/3/2,2/1/3/1}
				\draw[cyan] (2*\x,2*\y) -- (2*\z,2*\u);	
				\foreach \x/\y/\z/\u in {1/3/3/3}
				\draw[cyan] (2*\x,2*\y+0.4) -- (2*\z,2*\u+0.4);	
				\foreach \x/\y/\z/\u in {2/3/3/3}
				\draw[cyan] (2*\x,2*\y+0.8) -- (2*\z,2*\u+0.8);						
				
				\foreach \x/\y in {0/3,1/3,2/3,3/3,1/2,2/2,3/2,2/1,3/1,3/0}
				\draw[fill=black] (2*\x,\y*2) circle (3pt);
				
				\foreach \x/\y in {1/3,2/3,3/3,2/2,3/2,3/1}
				\draw[fill=blue] (2*\x,\y*2+0.4) circle (3pt);
				
				\foreach \x/\y in {2/3,3/3,3/2}
				\draw[fill=magenta] (2*\x,\y*2+0.8) circle (3pt);
				
				\foreach \x/\y in {3/3}
				\draw[fill=orange] (2*\x,\y*2+1.2) circle (3pt);					
				
		\end{tikzpicture}\end{center}\caption{\label{FigureWH00}Graphical presentation of domain and image of all possible weak homomorphisms $f:P_4\longrightarrow P_4\square P_5$ where $f(0)=(0,0)$.}
	\end{figure}

	Figure \ref{FigureWH00} shows all possible weak homomorphisms from  $P_4$ to $P_4\square P_5$ which map $0$ to $(0,0)$. The numbers on top are elements of domain set $V(P_4)$ and the tuples on the left are elements of image set $V(P_4\square P_5)$. The tuples with the same second elements are represented by the circle with the same color.
	
	We noted that normal black lines represent the increment of the first coordinate, dashed black lines represent the decrement of the first coordinate, normal magenta lines represent the increment of the second coordinate, magenta lines represent the decrement of the second coordinate and cyan lines represent no change in both coordinates.
	
	We now divide all the mappings in $\mathrm{Hom}^{00}(P_4,P_4 \square P_5)$ into groups according to the number of changes in the first coordinate $h$, and rewrite each path as 2 shorter paths.  The first path is formed by gray lines. On the other hand, the second path consists of cyan and magenta lines. In both paths, lines are arranged in sequential order.

	\begin{center}\begin{tabular}{|c|p{1.8in}|p{2.1in}|c|}
		\hline $h$ & $f\in \mathrm{WHom}^{00}(P_4,P_4\square P_5)$ with changes in the first coordinate $h$ times & \centering Paths represent each $f\in \mathrm{WHom}^{00}(P_4,P_4\square P_5)$ \qquad (Expanded Diagram) &  $P_{h+1}\ \&\ P_{4-h}$\\\hline
	\begin{minipage}[b]{0.3cm}$0$\vspace{3.5cm}\end{minipage} &	\begin{minipage}[b]{0.33\textwidth}
			\begin{tikzpicture}[scale=0.4]
				\draw (3.5,4+0.5) node[anchor=west] {$P_4$};
				\draw (-3,-0.5) node[anchor=north] {$P_4\square P_5$};
				
				\foreach \x in {0,...,3}
				{\draw[fill=black] (\x,4.5) circle (3pt);
					\draw (\x,4.2+0.5) node[anchor=south] {$\mathrm{\x}$};}
				\foreach \x in {0,...,2}
				{\draw (\x,4+0.5)-- ({\x+1},4+0.5);}
				
				\foreach \y in {0,...,3}
				\draw (-1,\y)-- (-5,\y);
				\foreach \x in {-1,...,-5}
				\draw (\x,0)-- (\x,3);				
				
				\foreach \y in {0,...,3}
				{\draw[fill=green] (-1,{-1*\y+3}) circle (3pt);
					\draw[fill=orange] (-2,{-1*\y+3}) circle (3pt);
					\draw[fill=magenta] (-3,{-1*\y+3}) circle (3pt);
					\draw[fill=blue] (-4,{-1*\y+3}) circle (3pt);
					\draw[fill=black] (-5,{-1*\y+3}) circle (3pt);}
				
				\foreach \x/\y/\z/\u in {0/3/1/3,1/3/2/3,2/3/3/3}
				\draw[magenta] (\x,\y) -- (\z,\u+0.25);								
				\foreach \x/\y/\z/\u in {1/3/2/3,2/3/3/3}
				\draw[magenta] (\x,\y+0.25) -- (\z,\u+0.5);									
				\foreach \x/\y/\z/\u in {2/3/3/3}
				\draw[magenta] (\x,\y+0.5) -- (\z,\u+0.75);	
				
				\foreach \x/\y/\z/\u in {1/3/2/3,2/3/3/3}
				\draw[magenta,dashed] (\x,\y+0.25) -- (\z,\u);						
				\foreach \x/\y/\z/\u in {2/3/3/3}
				\draw[magenta,dashed] (\x,\y+0.5) -- (\z,\u+0.25);	
				
				\foreach \x/\y/\z/\u in {0/3/3/3}
				\draw[cyan] (\x,\y) -- (\z,\u);	
				\foreach \x/\y/\z/\u in {1/3/3/3}
				\draw[cyan] (\x,\y+0.25) -- (\z,\u+0.25);	
				\foreach \x/\y/\z/\u in {2/3/3/3}
				\draw[cyan] (\x,\y+0.5) -- (\z,\u+0.5);						
				
				\foreach \x/\y in {0/3,1/3,2/3,3/3}
				\draw[fill=black] (\x,\y) circle (3pt);
				
				\foreach \x/\y in {1/3,2/3,3/3}
				\draw[fill=blue] (\x,\y+0.25) circle (3pt);
				
				\foreach \x/\y in {2/3,3/3}
				\draw[fill=magenta] (\x,\y+0.5) circle (3pt);
				
				\foreach \x/\y in {3/3}
				\draw[fill=orange] (\x,\y+0.75) circle (3pt);			
				
				\draw (0,-3.5) node[inner sep=1pt] {$\hfill$};
			\end{tikzpicture}
		\end{minipage}&
		\begin{minipage}[b]{0.33\textwidth}
			\begin{tikzpicture}[scale=0.4]				
				\def\l{2}
				\def\h{5}
				
				
				\draw[magenta] (0,3+2*\l) -- (1,3+0.25+2*\l);			
				\draw[magenta] (1,3+0.25+2*\l) -- (2,3+0.5+2*\l);		
				\draw[magenta] (2,3+0.5+2*\l) -- (3,3+0.75+2*\l);	
				\draw[fill=black] (0,3+2*\l) circle (3pt);
				\draw[fill=blue] (1,3+0.25+2*\l) circle (3pt);
				\draw[fill=magenta] (2,3+0.5+2*\l) circle (3pt);
				\draw[fill=orange] (3,3+0.75+2*\l) circle (3pt);
				
				\draw[magenta] (0+\h,3+2*\l) -- (1+\h,3+0.25+2*\l);			
				\draw[magenta] (1+\h,3+0.25+2*\l) -- (2+\h,3+0.5+2*\l);	
				\draw[cyan] (2+\h,3+0.5+2*\l) -- (3+\h,3+0.5+2*\l);	
				\draw[fill=black] (0+\h,3+2*\l) circle (3pt);
				\draw[fill=blue] (1+\h,3+0.25+2*\l) circle (3pt);
				\draw[fill=magenta] (2+\h,3+0.5+2*\l) circle (3pt);
				\draw[fill=magenta] (3+\h,3+0.5+2*\l) circle (3pt);
				
				\draw[magenta] (0+2*\h,3+2*\l) -- (1+2*\h,3+0.25+2*\l);			
				\draw[magenta] (1+2*\h,3+0.25+2*\l) -- (2+2*\h,3+0.5+2*\l);						
				\draw[magenta,dashed] (2+2*\h,3+0.5+2*\l) -- (3+2*\h,3+0.25+2*\l);											
				\draw[fill=black] (0+2*\h,3+2*\l) circle (3pt);
				\draw[fill=blue] (1+2*\h,3+0.25+2*\l) circle (3pt);
				\draw[fill=magenta] (2+2*\h,3+0.5+2*\l) circle (3pt);
				\draw[fill=blue] (3+2*\h,3+0.25+2*\l) circle (3pt);
				
				\draw[magenta] (0,3+\l) -- (1,3+0.25+\l);			
				\draw[cyan] (1,3+0.25+\l) -- (3,3+0.25+\l);	
				\draw[fill=black] (0,3+\l) circle (3pt);
				\draw[fill=blue] (1,3+0.25+\l) circle (3pt);
				\draw[fill=blue] (2,3+0.25+\l) circle (3pt);
				\draw[fill=blue] (3,3+0.25+\l) circle (3pt);
				
				\draw[magenta] (0+\h,3+\l) -- (1+\h,3+0.25+\l);			
				\draw[cyan] (1+\h,3+0.25+\l) -- (2+\h,3+0.25+\l);	
				\draw[magenta] (2+\h,3+0.25+\l) -- (3+\h,3+0.5+\l);	
				\draw[fill=black] (0+\h,3+\l) circle (3pt);
				\draw[fill=blue] (1+\h,3+0.25+\l) circle (3pt);
				\draw[fill=blue] (2+\h,3+0.25+\l) circle (3pt);
				\draw[fill=magenta] (3+\h,3+0.5+\l) circle (3pt);

				\draw[magenta] (0+2*\h,3+\l) -- (1+2*\h,3+0.25+\l);			
				\draw[cyan] (1+2*\h,3+0.25+\l) -- (2+2*\h,3+0.25+\l);	
				\draw[magenta,dashed] (2+2*\h,3+0.25+\l) -- (3+2*\h,3+\l);	
				\draw[fill=black] (0+2*\h,3+\l) circle (3pt);
				\draw[fill=blue] (1+2*\h,3+0.25+\l) circle (3pt);
				\draw[fill=blue] (2+2*\h,3+0.25+\l) circle (3pt);
				\draw[fill=black] (3+2*\h,3+\l) circle (3pt);

				\draw[magenta] (0,3) -- (1,3+0.25);			
				\draw[magenta,dashed] (1,3+0.25) -- (2,3);		
				\draw[cyan] (2,3) -- (3,3);	
				\draw[fill=black] (0,3) circle (3pt);
				\draw[fill=blue] (1,3+0.25) circle (3pt);
				\draw[fill=black] (2,3) circle (3pt);
				\draw[fill=black] (3,3) circle (3pt);
				
				\draw[magenta] (0+\h,3) -- (1+\h,3+0.25);			
				\draw[magenta,dashed] (1+\h,3+0.25) -- (2+\h,3);			
				\draw[magenta] (2+\h,3) -- (3+\h,3+0.25);	
				\draw[fill=black] (0+\h,3) circle (3pt);
				\draw[fill=blue] (1+\h,3+0.25) circle (3pt);
				\draw[fill=black] (2+\h,3) circle (3pt);
				\draw[fill=blue] (3+\h,3+0.25) circle (3pt);
				
				\draw[cyan] (0+2*\h,3) -- (1+2*\h,3);	
				\draw[magenta] (1+2*\h,3) -- (2+2*\h,3+0.25);	
				\draw[magenta] (2+2*\h,3+0.25) -- (3+2*\h,3+0.5);	
				\draw[fill=black] (0+2*\h,3) circle (3pt);
				\draw[fill=black] (1+2*\h,3) circle (3pt);
				\draw[fill=blue] (2+2*\h,3+0.25) circle (3pt);
				\draw[fill=magenta] (3+2*\h,3+0.5) circle (3pt);

				\draw[cyan] (0,3-\l) -- (1,3-\l);	
				\draw[magenta] (1,3-\l) -- (2,3+0.25-\l);	
				\draw[magenta,dashed] (2,3+0.25-\l) -- (3,3-\l);	
				\draw[fill=black] (0,3-\l) circle (3pt);
				\draw[fill=black] (1,3-\l) circle (3pt);
				\draw[fill=blue] (2,3+0.25-\l) circle (3pt);
				\draw[fill=black] (3,3-\l) circle (3pt);
				
				\draw[cyan] (0+\h,3-\l) -- (1+\h,3-\l);	
				\draw[magenta] (1+\h,3-\l) -- (2+\h,3+0.25-\l);	
				\draw[cyan] (2+\h,3+0.25-\l) -- (3+\h,3+0.25-\l);	
				\draw[fill=black] (0+\h,3-\l) circle (3pt);
				\draw[fill=black] (1+\h,3-\l) circle (3pt);
				\draw[fill=blue] (2+\h,3+0.25-\l) circle (3pt);
				\draw[fill=blue] (3+\h,3+0.25-\l) circle (3pt);
				
				\draw[cyan] (0+2*\h,3-\l) -- (2+2*\h,3-\l);						
				\draw[magenta] (2+2*\h,3-\l) -- (3+2*\h,3+0.25-\l);	
				\draw[fill=black] (0+2*\h,3-\l) circle (3pt);
				\draw[fill=black] (1+2*\h,3-\l) circle (3pt);
				\draw[fill=black] (2+2*\h,3-\l) circle (3pt);
				\draw[fill=blue] (3+2*\h,3+0.25-\l) circle (3pt);
				
				\draw[cyan] (0,3-2*\l) -- (3,3-2*\l);	
				\draw[fill=black] (0,3-2*\l) circle (3pt);
				\draw[fill=black] (1,3-2*\l) circle (3pt);
				\draw[fill=black] (2,3-2*\l) circle (3pt);
				\draw[fill=black] (3,3-2*\l) circle (3pt);
				
				\draw (0,4*\l+0.5) node[inner sep=1pt] {$\hfill$};
				\draw (0,-\l) node[inner sep=1pt] {$\hfill$};
			\end{tikzpicture}
		\end{minipage}&
		\begin{minipage}[b]{0.2\textwidth}
			\begin{tikzpicture}[scale=0.4]
				\def\l{2}
				\def\h{2}
				
				\draw[fill=gray] (0,0+5*\h) circle (3pt);
				\draw[magenta] (0+\l,0+5*\h)-- (1+\l,0+5*\h) -- (2+\l,0+5*\h)-- (3+\l,0+5*\h);
				\draw[fill=gray] (0+\l,0+5*\h) circle (3pt);
				\draw[fill=gray] (1+\l,0+5*\h) circle (3pt);
				\draw[fill=gray] (2+\l,0+5*\h) circle (3pt);		
				\draw[fill=gray] (3+\l,0+5*\h) circle (3pt);		
				\draw (1,-0.5+5*\h) node[inner sep=1pt,rectangle,fill=white] {, };
				
				\draw[fill=gray] (0,0+4*\h) circle (3pt);
				\draw[magenta] (0+\l,0+4*\h)-- (1+\l,0+4*\h) -- (2+\l,0+4*\h);
				\draw[cyan] (2+\l,0+4*\h)-- (3+\l,0+4*\h);
				\draw[fill=gray] (0+\l,0+4*\h) circle (3pt);
				\draw[fill=gray] (1+\l,0+4*\h) circle (3pt);
				\draw[fill=gray] (2+\l,0+4*\h) circle (3pt);		
				\draw[fill=gray] (3+\l,0+4*\h) circle (3pt);		
				\draw (1,-0.5+4*\h) node[inner sep=1pt,rectangle,fill=white] {, };
				
				\draw[fill=gray] (0,0+3*\h) circle (3pt);
				\draw[magenta] (0+\l,0+3*\h)-- (1+\l,0+3*\h) -- (2+\l,0+3*\h);
				\draw[dashed, red] (2+\l,0+3*\h) -- (3+\l,0+3*\h);
				\draw[fill=gray] (0+\l,0+3*\h) circle (3pt);
				\draw[fill=gray] (1+\l,0+3*\h) circle (3pt);
				\draw[fill=gray] (2+\l,0+3*\h) circle (3pt);		
				\draw[fill=gray] (3+\l,0+3*\h) circle (3pt);		
				\draw (1,-0.5+3*\h) node[inner sep=1pt,rectangle,fill=white] {, };

				\draw (1,0+2*\h) node[inner sep=1pt,rectangle,fill=white] {$\vdots$};
				
				\draw[fill=gray] (0,0+\h) circle (3pt);
				\draw[cyan] (0+\l,0+\h)-- (1+\l,0+\h) -- (2+\l,0+\h)-- (3+\l,0+\h);
				\draw[fill=gray] (0+\l,0+\h) circle (3pt);
				\draw[fill=gray] (1+\l,0+\h) circle (3pt);
				\draw[fill=gray] (2+\l,0+\h) circle (3pt);		
				\draw[fill=gray] (3+\l,0+\h) circle (3pt);		
				\draw (1,-0.5+\h) node[inner sep=1pt,rectangle,fill=white] {, };
				
				\draw (0,0.75) node[inner sep=1pt] {$\hfill$};
				
			\end{tikzpicture}
		\end{minipage}
			\\\hline \begin{minipage}[b]{0.3cm}$1$\vspace{3.7cm}\end{minipage} &
			\begin{minipage}[b]{0.33\textwidth}
				\begin{tikzpicture}[scale=0.4]
					\draw (3.5,4+0.5) node[anchor=west] {$P_4$};
					\draw (-3,-0.5) node[anchor=north] {$P_4\square P_5$};
					
					\foreach \x in {0,...,3}
					{\draw[fill=black] (\x,4.5) circle (3pt);
						\draw (\x,4.2+0.5) node[anchor=south] {$\mathrm{\x}$};}
					\foreach \x in {0,...,2}
					{\draw (\x,4+0.5)-- ({\x+1},4+0.5);}
					
					\foreach \y in {0,...,3}
					\draw (-1,\y)-- (-5,\y);
					\foreach \x in {-1,...,-5}
					\draw (\x,0)-- (\x,3);				
					
					\foreach \y in {0,...,3}
					{\draw[fill=green] (-1,{-1*\y+3}) circle (3pt);
						\draw[fill=orange] (-2,{-1*\y+3}) circle (3pt);
						\draw[fill=magenta] (-3,{-1*\y+3}) circle (3pt);
						\draw[fill=blue] (-4,{-1*\y+3}) circle (3pt);
						\draw[fill=black] (-5,{-1*\y+3}) circle (3pt);}
					
					\foreach \x/\y/\z/\u in {0/3/1/2,1/3/2/2,2/3/3/2}
					\draw[gray] (\x,\y) -- (\z,\u);			
					\foreach \x/\y/\z/\u in {1/3/2/2,2/3/3/2}
					\draw[gray] (\x,\y+0.25) -- (\z,\u+0.25);			
					\foreach \x/\y/\z/\u in {2/3/3/2}
					\draw[gray] (\x,\y+0.5) -- (\z,\u+0.5);	
					
					\draw[magenta] (0,3) -- (1,3+0.25);			
					\draw[magenta] (1,3) -- (2,3+0.25);	
					\draw[magenta] (1,2) -- (2,2+0.25);		
					\draw[magenta] (2,2) -- (3,2+0.25);					
					
					\draw[magenta] (1,3+0.25) -- (2,3+0.5);		
					\draw[magenta] (2,2+0.25) -- (3,2+0.5);							
					
					\draw[magenta,dashed] (1,3+0.25) -- (2,3);	
					\draw[magenta,dashed] (2,2+0.25) -- (3,2);					
					
					\draw[cyan] (0,3) -- (1,3);	
					\draw[cyan] (1,3) -- (2,3);	
					\draw[cyan] (1,2) -- (2,2);	
					\draw[cyan] (2,2) -- (3,2);	
					
					\draw[cyan] (1,3+0.25) -- (2,3+0.25);	
					\draw[cyan] (2,2+0.25) -- (3,2+0.25);	
					
					\foreach \x/\y in {0/3,1/3,2/3,1/2,2/2,3/2}
					\draw[fill=black] (\x,\y) circle (3pt);
					
					\foreach \x/\y in {1/3,2/3,2/2,3/2}
					\draw[fill=blue] (\x,\y+0.25) circle (3pt);
					
					\foreach \x/\y in {2/3,3/2}
					\draw[fill=magenta] (\x,\y+0.5) circle (3pt);
					
					\draw (0,-3.5) node[inner sep=1pt] {$\hfill$};
					
				\end{tikzpicture}
			\end{minipage}&
			\begin{minipage}[b]{0.45\textwidth}
				\begin{tikzpicture}[scale=0.4]
					\def\l{5}
					\def\h{2}
					
					\draw[magenta] (0,3+2*\h) -- (1,3+0.25+2*\h);	
					\draw[magenta] (1,3+0.25+2*\h) -- (2,3+0.5+2*\h);		
					\draw[gray] (2,3+0.5+2*\h) -- (3,2+0.5+2*\h);	
					\draw[fill=black] (0,3+2*\h) circle (3pt);
					\draw[fill=blue] (1,3+0.25+2*\h) circle (3pt);
					\draw[fill=magenta] (2,3+0.5+2*\h) circle (3pt);
					\draw[fill=magenta] (3,2+0.5+2*\h) circle (3pt);
					
					\draw[gray] (0+\l,3+2*\h) -- (1+\l,2+2*\h);			
					\draw[magenta] (1+\l,2+2*\h) -- (2+\l,2+0.25+2*\h);		
					\draw[magenta] (2+\l,2+0.25+2*\h) -- (3+\l,2+0.5+2*\h);	
					\draw[fill=black] (0+\l,3+2*\h) circle (3pt);
					\draw[fill=black] (1+\l,2+2*\h) circle (3pt);
					\draw[fill=blue] (2+\l,2+0.25+2*\h) circle (3pt);
					\draw[fill=magenta] (3+\l,2+0.5+2*\h) circle (3pt);					
					
					\draw[magenta] (0+2*\l,3+2*\h) -- (1+2*\l,3+0.25+2*\h);	
					\draw[gray] (1+2*\l,3+0.25+2*\h) -- (2+2*\l,2+0.25+2*\h);	
					\draw[magenta] (2+2*\l,2+0.25+2*\h) -- (3+2*\l,2+0.5+2*\h);	
					\draw[fill=black] (0+2*\l,3+2*\h) circle (3pt);
					\draw[fill=blue] (1+2*\l,3+0.25+2*\h) circle (3pt);
					\draw[fill=blue] (2+2*\l,2+0.25+2*\h) circle (3pt);
					\draw[fill=magenta] (3+2*\l,2+0.5+2*\h) circle (3pt);
					
					\draw[magenta] (0,3+\h) -- (1,3+0.25+\h);	
					\draw[magenta,dashed] (1,3+0.25+\h) -- (2,3+\h);	
					\draw[gray] (2,3+\h) -- (3,2+\h);	
					\draw[fill=black] (0,3+\h) circle (3pt);
					\draw[fill=blue] (1,3+0.25+\h) circle (3pt);
					\draw[fill=black] (2,3+\h) circle (3pt);
					\draw[fill=black] (3,2+\h) circle (3pt);
					
					\draw[magenta] (0+\l,3+\h) -- (1+\l,3+0.25+\h);	
					\draw[gray] (1+\l,3+0.25+\h) -- (2+\l,2+0.25+\h);	
					\draw[magenta,dashed] (2+\l,2+0.25+\h) -- (3+\l,2+\h);
					\draw[fill=black] (0+\l,3+\h) circle (3pt);
					\draw[fill=blue] (1+\l,3+0.25+\h) circle (3pt);
					\draw[fill=blue] (2+\l,2+0.25+\h) circle (3pt);
					\draw[fill=black] (3+\l,2+\h) circle (3pt);
					
					\draw[gray] (0+2*\l,3+\h) -- (1+2*\l,2+\h);			
					\draw[magenta] (1+2*\l,2+\h) -- (2+2*\l,2+0.25+\h);		
					\draw[magenta,dashed] (2+2*\l,2+0.25+\h) -- (3+2*\l,2+\h);
					\draw[fill=black] (0+2*\l,3+\h) circle (3pt);
					\draw[fill=black] (1+2*\l,2+\h) circle (3pt);
					\draw[fill=blue] (2+2*\l,2+0.25+\h) circle (3pt);
					\draw[fill=black] (3+2*\l,2+\h) circle (3pt);
					
					\draw[magenta] (0,3) -- (1,3+0.25);	
					\draw[cyan] (1,3+0.25) -- (2,3+0.25);	
					\draw[gray] (2,3+0.25) -- (3,2+0.25);
					\draw[fill=black] (0,3) circle (3pt);
					\draw[fill=blue] (1,3+0.25) circle (3pt);
					\draw[fill=blue] (2,3+0.25) circle (3pt);
					\draw[fill=blue] (3,2+0.25) circle (3pt);
					
					\draw[magenta] (0+\l,3) -- (1+\l,3+0.25);	
					\draw[gray] (1+\l,3+0.25) -- (2+\l,2+0.25);	
					\draw[cyan] (2+\l,2+0.25) -- (3+\l,2+0.25);	
					\draw[fill=black] (0+\l,3) circle (3pt);
					\draw[fill=blue] (1+\l,3+0.25) circle (3pt);
					\draw[fill=blue] (2+\l,2+0.25) circle (3pt);
					\draw[fill=blue] (3+\l,2+0.25) circle (3pt);
					
					\draw[gray] (0+2*\l,3) -- (1+2*\l,2);	
					\draw[magenta] (1+2*\l,2) -- (2+2*\l,2+0.25);	
					\draw[cyan] (2+2*\l,2+0.25) -- (3+2*\l,2+0.25);	
					\draw[fill=black] (0+2*\l,3) circle (3pt);
					\draw[fill=black] (1+2*\l,2) circle (3pt);
					\draw[fill=blue] (2+2*\l,2+0.25) circle (3pt);
					\draw[fill=blue] (3+2*\l,2+0.25) circle (3pt);
					
					\draw[cyan] (0,3-\h) -- (1,3-\h);	
					\draw[gray] (1,3-\h) -- (2,2-\h);	
					\draw[magenta] (2,2-\h) -- (3,2+0.25-\h);	
					\draw[fill=black] (0,3-\h) circle (3pt);
					\draw[fill=black] (1,3-\h) circle (3pt);
					\draw[fill=black] (2,2-\h) circle (3pt);
					\draw[fill=blue] (3,2+0.25-\h) circle (3pt);
					
					\draw[cyan] (0+\l,3-\h) -- (1+\l,3-\h);	
					\draw[magenta] (1+\l,3-\h) -- (2+\l,3+0.25-\h);	
					\draw[gray] (2+\l,3+0.25-\h) -- (3+\l,2+0.25-\h);
					\draw[fill=black] (0+\l,3-\h) circle (3pt);
					\draw[fill=black] (1+\l,3-\h) circle (3pt);
					\draw[fill=blue] (2+\l,3+0.25-\h) circle (3pt);
					\draw[fill=blue] (3+\l,2+0.25-\h) circle (3pt);
					
					\draw[gray] (0+2*\l,3-\h) -- (1+2*\l,2-\h);			
					\draw[cyan] (1+2*\l,2-\h) -- (2+2*\l,2-\h);	
					\draw[magenta] (2+2*\l,2-\h) -- (3+2*\l,2+0.25-\h);	
					\draw[fill=black] (0+2*\l,3-\h) circle (3pt);
					\draw[fill=black] (1+2*\l,2-\h) circle (3pt);
					\draw[fill=black] (2+2*\l,2-\h) circle (3pt);
					\draw[fill=blue] (3+2*\l,2+0.25-\h) circle (3pt);
					
					\draw[cyan] (0,3-2*\h) -- (1,3-2*\h);	
					\draw[cyan] (1,3-2*\h) -- (2,3-2*\h);	
					\draw[gray] (2,3-2*\h) -- (3,2-2*\h);	
					\draw[fill=black] (0,3-2*\h) circle (3pt);
					\draw[fill=black] (1,3-2*\h) circle (3pt);
					\draw[fill=black] (2,3-2*\h) circle (3pt);
					\draw[fill=black] (3,2-2*\h) circle (3pt);
					
					\draw[gray] (0+\l,3-2*\h) -- (1+\l,2-2*\h);			
					\draw[cyan] (1+\l,2-2*\h) -- (2+\l,2-2*\h);	
					\draw[cyan] (2+\l,2-2*\h) -- (3+\l,2-2*\h);	
					\draw[fill=black] (0+\l,3-2*\h) circle (3pt);
					\draw[fill=black] (1+\l,2-2*\h) circle (3pt);
					\draw[fill=black] (2+\l,2-2*\h) circle (3pt);
					\draw[fill=black] (3+\l,2-2*\h) circle (3pt);
					
					\draw[cyan] (0+2*\l,3-2*\h) -- (1+2*\l,3-2*\h);	
					\draw[gray] (1+2*\l,3-2*\h) -- (2+2*\l,2-2*\h);	
					\draw[cyan] (2+2*\l,2-2*\h) -- (3+2*\l,2-2*\h);	
					\draw[fill=black] (0+2*\l,3-2*\h) circle (3pt);
					\draw[fill=black] (1+2*\l,3-2*\h) circle (3pt);
					\draw[fill=black] (2+2*\l,2-2*\h) circle (3pt);
					\draw[fill=black] (3+2*\l,2-2*\h) circle (3pt);
					
					\draw (0,4*\h+0.5) node[inner sep=1pt] {$\hfill$};
					\draw (0,-\h-0.5) node[inner sep=1pt] {$\hfill$};
				\end{tikzpicture}
			\end{minipage}&
			\begin{minipage}[b]{0.2\textwidth}
				\begin{tikzpicture}[scale=0.4]				
					\def\l{3}
					\def\h{2}
					\draw[gray] (0,0+2*\h)-- (1,0+2*\h);
					\draw[fill=gray] (0,0+2*\h) circle (3pt);
					\draw[fill=gray] (1,0+2*\h) circle (3pt);		
					\draw[magenta] (0+\l,0+2*\h)-- (1+\l,0+2*\h) -- (2+\l,0+2*\h);
					\draw[fill=gray] (0+\l,0+2*\h) circle (3pt);
					\draw[fill=gray] (1+\l,0+2*\h) circle (3pt);
					\draw[fill=gray] (2+\l,0+2*\h) circle (3pt);		
					\draw (2,-0.5+2*\h) node[inner sep=1pt,rectangle,fill=white] {, };
					
					\draw[gray] (0,0+\h)-- (1,0+\h);
					\draw[fill=gray] (0,0+\h) circle (3pt);
					\draw[fill=gray] (1,0+\h) circle (3pt);		
					\draw[magenta] (0+\l,0+\h)-- (1+\l,0+\h);
					\draw[magenta,dashed] (1+\l,0+\h) -- (2+\l,0+\h);
					\draw[fill=gray] (0+\l,0+\h) circle (3pt);
					\draw[fill=gray] (1+\l,0+\h) circle (3pt);
					\draw[fill=gray] (2+\l,0+\h) circle (3pt);		
					\draw (2,-0.5+\h) node[inner sep=1pt,rectangle,fill=white] {, };

					\draw[gray] (0,0)-- (1,0);
					\draw[fill=gray] (0,0) circle (3pt);
					\draw[fill=gray] (1,0) circle (3pt);		
					\draw[magenta] (0+\l,0)-- (1+\l,0);
					\draw[cyan] (1+\l,0) -- (2+\l,0);
					\draw[fill=gray] (0+\l,0) circle (3pt);
					\draw[fill=gray] (1+\l,0) circle (3pt);
					\draw[fill=gray] (2+\l,0) circle (3pt);		
					\draw (2,-0.5) node[inner sep=1pt,rectangle,fill=white] {, };
					\draw (0,-1) node[inner sep=1pt,fill=white] {$\hfill$};
					
					\draw[gray] (0,0-\h)-- (1,0-\h);
					\draw[fill=gray] (0,0-\h) circle (3pt);
					\draw[fill=gray] (1,0-\h) circle (3pt);		
					\draw[cyan] (0+\l,0-\h)-- (1+\l,0-\h) ;
					\draw[magenta] (1+\l,0-\h) -- (2+\l,0-\h);
					\draw[fill=gray] (0+\l,0-\h) circle (3pt);
					\draw[fill=gray] (1+\l,0-\h) circle (3pt);
					\draw[fill=gray] (2+\l,-+\h) circle (3pt);		
					\draw (2,-0.5-\h) node[inner sep=1pt,rectangle,fill=white] {, };
					
					\draw[gray] (0,0-2*\h)-- (1,0-2*\h);
					\draw[fill=gray] (0,0-2*\h) circle (3pt);
					\draw[fill=gray] (1,0-2*\h) circle (3pt);		
					\draw[cyan] (0+\l,0-2*\h)-- (1+\l,0-2*\h) -- (2+\l,0-2*\h);
					\draw[fill=gray] (0+\l,0-2*\h) circle (3pt);
					\draw[fill=gray] (1+\l,0-2*\h) circle (3pt);
					\draw[fill=gray] (2+\l,0-2*\h) circle (3pt);		
					\draw (2,-0.5-2*\h) node[inner sep=1pt,rectangle,fill=white] {, };
					
					\draw (0,-3*\h+1) node[inner sep=1pt] {$\hfill$};

				\end{tikzpicture}
			\end{minipage}
			\\
			\hline
			
			\begin{minipage}[b]{0.3cm}$2$\vspace{3.2cm}\end{minipage} &	\begin{minipage}[b]{0.33\textwidth}
				\begin{tikzpicture}[scale=0.4]
					\draw (3.5,4+0.5) node[anchor=west] {$P_4$};
					\draw (-3,-0.5) node[anchor=north] {$P_4\square P_5$};
					
					\foreach \x in {0,...,3}
					{\draw[fill=black] (\x,4.5) circle (3pt);
						\draw (\x,4.2+0.5) node[anchor=south] {$\mathrm{\x}$};}
					\foreach \x in {0,...,2}
					{\draw (\x,4+0.5)-- ({\x+1},4+0.5);}
					
					\foreach \y in {0,...,3}
					\draw (-1,\y)-- (-5,\y);
					\foreach \x in {-1,...,-5}
					\draw (\x,0)-- (\x,3);				
					
					\foreach \y in {0,...,3}
					{\draw[fill=green] (-1,{-1*\y+3}) circle (3pt);
						\draw[fill=orange] (-2,{-1*\y+3}) circle (3pt);
						\draw[fill=magenta] (-3,{-1*\y+3}) circle (3pt);
						\draw[fill=blue] (-4,{-1*\y+3}) circle (3pt);
						\draw[fill=black] (-5,{-1*\y+3}) circle (3pt);}

					\foreach \x/\y/\z/\u in {0/3/2/1,1/3/3/1}
					\draw[gray] (\x,\y) -- (\z,\u);			
					\foreach \x/\y/\z/\u in {1/3/3/1}
					\draw[gray] (\x,\y+0.25) -- (\z,\u+0.25);

					\foreach \x/\y/\z/\u in {1/2/2/3,2/2/3/3}
					\draw[gray,dashed] (\x,\y) -- (\z,\u);			
					\foreach \x/\y/\z/\u in {2/2/3/3}
					\draw[gray,dashed] (\x,\y+0.25) -- (\z,\u+0.25);			
					
					\foreach \x/\y/\z/\u in {0/3/1/3,2/3/3/3,1/2/2/2,2/1/3/1}
					\draw[magenta] (\x,\y) -- (\z,\u+0.25);

					\foreach \x/\y/\z/\u in {0/3/1/3,1/2/2/2,2/1/3/1,2/3/3/3}
					\draw[cyan] (\x,\y) -- (\z,\u);

					\foreach \x/\y in {0/3,1/3,2/3,3/3,1/2,2/2,2/1,3/1}
					\draw[fill=black] (\x,\y) circle (3pt);
					
					\foreach \x/\y in {1/3,3/3,2/2,3/1}
					\draw[fill=blue] (\x,\y+0.25) circle (3pt);
					
					\draw (0,-3) node[inner sep=1pt] {$\hfill$};

				\end{tikzpicture}
			\end{minipage}&
			\begin{minipage}[b]{0.45\textwidth}
				\begin{tikzpicture}[scale=0.4]
					\def\l{5}
					\def\h{2}
					
					\draw[magenta] (0,3+2*\h) -- (1,3+0.25+2*\h);		
					\draw[gray] (1,3+0.25+2*\h) -- (2,2+0.25+2*\h);
					\draw[gray,dashed] (2,2+0.25+2*\h) -- (3,3+0.25+2*\h);
					\draw[fill=black] (0,3+2*\h) circle (3pt);
					\draw[fill=blue] (1,3+0.25+2*\h) circle (3pt);
					\draw[fill=black] (2,2+0.25+2*\h) circle (3pt);
					\draw[fill=black] (3,3+0.25+2*\h) circle (3pt);
					
					\draw[gray] (0+\l,3+2*\h) -- (1+\l,2+2*\h);	
					\draw[magenta] (1+\l,2+2*\h) -- (2+\l,2+0.25+2*\h);
					\draw[gray,dashed] (2+\l,2+0.25+2*\h) -- (3+\l,3+0.25+2*\h);
					\draw[fill=black] (0+\l,3+2*\h) circle (3pt);
					\draw[fill=black] (1+\l,2+2*\h)  circle (3pt);
					\draw[fill=blue] (2+\l,2+0.25+2*\h) circle (3pt);
					\draw[fill=black] (3+\l,3+0.25+2*\h)  circle (3pt);
					
					\draw[gray] (0+2*\l,3+2*\h) -- (1+2*\l,2+2*\h);	
					\draw[gray,dashed] (1+2*\l,2+2*\h) -- (2+2*\l,3+2*\h);
					\draw[magenta] (2+2*\l,3+2*\h) -- (3+2*\l,3+0.25+2*\h);
					\draw[fill=black] (0+2*\l,3+2*\h) circle (3pt);
					\draw[fill=black] (1+2*\l,2+2*\h) circle (3pt);
					\draw[fill=black] (2+2*\l,3+2*\h) circle (3pt);
					\draw[fill=blue]  (3+2*\l,3+0.25+2*\h) circle (3pt);
					
					\def\h{2.5}
					\draw[magenta] (0,3+\h) -- (1,3+0.25+\h);		
					\draw[gray] (1,3+0.25+\h) -- (2,2+0.25+\h);
					\draw[gray] (2,2+0.25+\h) -- (3,1+0.25+\h);
					\draw[fill=black] (0,3+\h) circle (3pt);
					\draw[fill=blue] (1,3+0.25+\h) circle (3pt);
					\draw[fill=black]  (2,2+0.25+\h) circle (3pt);
					\draw[fill=black] (3,1+0.25+\h) circle (3pt);
					
					\draw[gray] (0+\l,3+\h) -- (2+\l,1+\h);	
					\draw[magenta] (2+\l,1+\h) -- (3+\l,1+0.25+\h);	
					\draw[fill=black] (0+\l,3+\h) circle (3pt);
					\draw[fill=black] (1+\l,2+\h) circle (3pt);
					\draw[fill=black] (2+\l,1+\h) circle (3pt);
					\draw[fill=blue] (3+\l,1+0.25+\h) circle (3pt);
					
					\draw[gray] (0+2*\l,3+\h) -- (1+2*\l,2+\h);	
					\draw[magenta] (1+2*\l,2+\h) -- (2+2*\l,2+0.25+\h);
					\draw[gray] (2+2*\l,2+0.25+\h) -- (3+2*\l,1+0.25+\h);
					\draw[fill=black](0+2*\l,3+\h)  circle (3pt);
					\draw[fill=black] (1+2*\l,2+\h) circle (3pt);
					\draw[fill=blue] (2+2*\l,2+0.25+\h) circle (3pt);
					\draw[fill=black]  (3+2*\l,1+0.25+\h) circle (3pt);

					\draw[gray] (0,3) -- (1,2);	
					\draw[gray,dashed] (1,2) -- (2,3);
					\draw[cyan] (2,3) -- (3,3);
					\draw[fill=black] (0,3) circle (3pt);
					\draw[fill=black] (1,2) circle (3pt);
					\draw[fill=black]  (2,3) circle (3pt);
					\draw[fill=black]  (3,3) circle (3pt);
					
					\draw[cyan] (0+\l,3) -- (1+\l,3);			
					\draw[gray] (1+\l,3) -- (2+\l,2);
					\draw[gray,dashed] (2+\l,2) -- (3+\l,3);	
					\draw[fill=black] (0+\l,3) circle (3pt);
					\draw[fill=black] (1+\l,3) circle (3pt);
					\draw[fill=black] (2+\l,2) circle (3pt);
					\draw[fill=black] (3+\l,3) circle (3pt);

					\draw[gray] (0+2*\l,3) -- (1+2*\l,2);	
					\draw[cyan] (1+2*\l,2) -- (2+2*\l,2);	
					\draw[gray,dashed] (2+2*\l,2) -- (3+2*\l,3);	
					\draw[fill=black](0+2*\l,3) circle (3pt);
					\draw[fill=black] (1+2*\l,2) circle (3pt);
					\draw[fill=black] (2+2*\l,2) circle (3pt);
					\draw[fill=black] (3+2*\l,3) circle (3pt);
					
					\def\h{2}
					\draw[gray] (0,3-\h) -- (1,2-\h);	
					\draw[cyan] (1,2-\h) -- (2,2-\h);	
					\draw[gray] (2,2-\h) -- (3,1-\h);	
					\draw[fill=black] (0,3-\h) circle (3pt);
					\draw[fill=black](1,2-\h)  circle (3pt);
					\draw[fill=black] (2,2-\h) circle (3pt);
					\draw[fill=black](3,1-\h)  circle (3pt);

					\draw[gray] (0+\l,3-\h) -- (2+\l,1-\h);	
					\draw[cyan] (2+\l,1-\h) -- (3+\l,1-\h);	
					\draw[fill=black] (0+\l,3-\h) circle (3pt);
					\draw[fill=black](1+\l,2-\h)  circle (3pt);
					\draw[fill=black] (2+\l,1-\h) circle (3pt);
					\draw[fill=black]  (3+\l,1-\h) circle (3pt);
					
					\draw[cyan] (0+2*\l,3-\h) -- (1+2*\l,3-\h);			
					\draw[gray] (1+2*\l,3-\h) -- (2+2*\l,2-\h);
					\draw[gray] (2+2*\l,2-\h) -- (3+2*\l,1-\h);
					\draw[fill=black] (0+2*\l,3-\h)  circle (3pt);
					\draw[fill=black] (1+2*\l,3-\h) circle (3pt);
					\draw[fill=black] (2+2*\l,2-\h) circle (3pt);
					\draw[fill=black] (3+2*\l,1-\h) circle (3pt);

					\draw (0,4*\h) node[inner sep=1pt,fill=white] {$\hfill$};
					\draw (0,-\h) node[inner sep=1pt,fill=white] {$\hfill$};
					
				\end{tikzpicture}
			\end{minipage}&
			\begin{minipage}[b]{0.2\textwidth}
				\begin{tikzpicture}[scale=0.4]				
					\def\l{4}
					\def\h{2.5}
					\draw[gray] (0,0+2*\h)-- (1,0+2*\h);
					\draw[gray,dashed] (1,0+2*\h)-- (2,0+2*\h);
					\draw[fill=gray] (0,0+2*\h) circle (3pt);
					\draw[fill=gray] (1,0+2*\h) circle (3pt);		
					\draw[fill=gray] (2,0+2*\h) circle (3pt);
					\draw[magenta] (0+\l,0+2*\h)-- (1+\l,0+2*\h) ;
					\draw[fill=gray] (0+\l,0+2*\h) circle (3pt);
					\draw[fill=gray] (1+\l,0+2*\h) circle (3pt);
					
					\draw (3,-0.5+2*\h) node[inner sep=1pt,rectangle,fill=white] {, };
					
					\draw[gray] (0,0+\h)-- (1,0+\h)-- (2,0+\h);
					\draw[fill=gray] (0,0+\h) circle (3pt);
					\draw[fill=gray] (1,0+\h) circle (3pt);		
					\draw[fill=gray] (2,0+\h) circle (3pt);
					\draw[magenta] (0+\l,0+\h)-- (1+\l,0+\h) ;
					\draw[fill=gray] (0+\l,0+\h) circle (3pt);
					\draw[fill=gray] (1+\l,0+\h) circle (3pt);
					
					\draw (3,-0.5+\h) node[inner sep=1pt,rectangle,fill=white] {, };
					
					\draw[gray] (0,0)-- (1,0);
					\draw[gray,dashed] (1,0)-- (2,0);
					\draw[fill=gray] (0,0) circle (3pt);
					\draw[fill=gray] (1,0) circle (3pt);
					\draw[fill=gray] (2,0) circle (3pt);			
					\draw[cyan] (0+\l,0)-- (1+\l,0);
					\draw[fill=gray] (0+\l,0) circle (3pt);
					\draw[fill=gray] (1+\l,0) circle (3pt);
					
					\draw (3,-0.5) node[inner sep=1pt,rectangle,fill=white] {, };
					
					\draw[gray] (0,0-\h)-- (1,0-\h);
					\draw[gray] (1,0-\h)-- (2,0-\h);
					\draw[fill=gray] (0,0-\h) circle (3pt);
					\draw[fill=gray] (1,0-\h) circle (3pt);
					\draw[fill=gray] (2,0-\h) circle (3pt);			
					\draw[cyan] (0+\l,0-\h)-- (1+\l,0-\h);
					\draw[fill=gray] (0+\l,0-\h) circle (3pt);
					\draw[fill=gray] (1+\l,0-\h) circle (3pt);
					
					\draw (3,-0.5-\h) node[inner sep=1pt,rectangle,fill=white] {, };
					
					\draw (0,-3.5) node[inner sep=1pt] {$\hfill$};
				\end{tikzpicture}
			\end{minipage}
			\\\hline
			\begin{minipage}[b]{0.3cm}$3$\vspace{2.5cm}\end{minipage} &
			\begin{minipage}[b]{0.3\textwidth}
				\begin{tikzpicture}[scale=0.4]
					\draw (3.5,4+0.5) node[anchor=west] {$P_4$};
					\draw (-3,-0.5) node[anchor=north] {$P_4\square P_5$};
					
					\foreach \x in {0,...,3}
					{\draw[fill=black] (\x,4.5) circle (3pt);
						\draw (\x,4.2+0.5) node[anchor=south] {$\mathrm{\x}$};}
					\foreach \x in {0,...,2}
					{\draw (\x,4+0.5)-- ({\x+1},4+0.5);}
					
					\foreach \y in {0,...,3}
					\draw (-1,\y)-- (-5,\y);
					\foreach \x in {-1,...,-5}
					\draw (\x,0)-- (\x,3);				
					
					\foreach \y in {0,...,3}
					{\draw[fill=green] (-1,{-1*\y+3}) circle (3pt);
						\draw[fill=orange] (-2,{-1*\y+3}) circle (3pt);
						\draw[fill=magenta] (-3,{-1*\y+3}) circle (3pt);
						\draw[fill=blue] (-4,{-1*\y+3}) circle (3pt);
						\draw[fill=black] (-5,{-1*\y+3}) circle (3pt);}
					
					\foreach \x/\y/\z/\u in {0/3/3/0,2/3/3/2}
					\draw[gray] (\x,\y) -- (\z,\u);

					\foreach \x/\y/\z/\u in {1/2/2/3,2/1/3/2}
					\draw[gray,dashed] (\x,\y) -- (\z,\u);

					\foreach \x/\y in {0/3,2/3,1/2,3/2,2/1,3/0}
					\draw[fill=black] (\x,\y) circle (3pt);

				\end{tikzpicture}
			\end{minipage}&
			\begin{minipage}[b]{0.45\textwidth}
				\begin{tikzpicture}[scale=0.4]
					
					\def\h{5}
					\draw[gray] (0,3)-- (1,2) -- (2,1) -- (3,0);
					\draw[fill=black] (0,3) circle (3pt);
					\draw[fill=black] (1,2) circle (3pt);
					\draw[fill=black] (2,1) circle (3pt);
					\draw[fill=black] (3,0) circle (3pt);
					
					\draw[gray] (0+\h,3)-- (1+\h,2);
					\draw[gray,dashed] (1+\h,2)-- (2+\h,3);
					\draw[gray] (2+\h,3)-- (3+\h,2);
					\draw[fill=black] (0+\h,3) circle (3pt);
					\draw[fill=black] (1+\h,2) circle (3pt);
					\draw[fill=black] (2+\h,3) circle (3pt);
					\draw[fill=black] (3+\h,2) circle (3pt);

					\draw[gray] (0+2*\h,3)-- (1+2*\h,2) -- (2+2*\h,1);
					\draw[gray,dashed] (2+2*\h,1)-- (3+2*\h,2);
					\draw[fill=black] (0+2*\h,3) circle (3pt);
					\draw[fill=black] (1+2*\h,2) circle (3pt);
					\draw[fill=black] (2+2*\h,1) circle (3pt);
					\draw[fill=black] (3+2*\h,2) circle (3pt);
					\draw (0,-1) node[inner sep=1pt,fill=white] {$\hfill$};
				\end{tikzpicture}
			\end{minipage}&
			\begin{minipage}[b]{0.2\textwidth}
				\begin{tikzpicture}[scale=0.4]				
					\def\l{5}
					\def\h{2}
					\draw[gray] (0,0+2*\h)-- (1,0+2*\h)-- (2,0+2*\h)-- (3,0+2*\h);
					\draw[fill=gray] (0,0+2*\h) circle (3pt);
					\draw[fill=gray] (1,0+2*\h) circle (3pt);		
					\draw[fill=gray] (2,0+2*\h) circle (3pt);
					\draw[fill=gray] (3,0+2*\h) circle (3pt);
					\draw[fill=gray] (0+\l,0+2*\h) circle (3pt);
					\draw (4,-0.5+2*\h) node[inner sep=1pt,rectangle,fill=white] {, };
					
					\draw[gray] (0,0+\h)-- (1,0+\h);
					\draw[gray]  (2,0+\h)-- (3,0+\h);
					\draw[gray,dashed]  (1,0+\h)-- (2,0+\h);
					\draw[fill=gray] (0,0+\h) circle (3pt);
					\draw[fill=gray] (1,0+\h) circle (3pt);		
					\draw[fill=gray] (2,0+\h) circle (3pt);
					\draw[fill=gray] (3,0+\h) circle (3pt);
					\draw[fill=gray] (0+\l,0+\h) circle (3pt);
					\draw (4,-0.5+\h) node[inner sep=1pt,rectangle,fill=white] {, };
					
					\draw[gray] (0,0)-- (1,0)-- (2,0);
					\draw[gray,dashed] (2,0)-- (3,0);
					\draw[fill=gray] (0,0) circle (3pt);
					\draw[fill=gray] (1,0) circle (3pt);		
					\draw[fill=gray] (2,0) circle (3pt);
					\draw[fill=gray] (3,0) circle (3pt);
					\draw[fill=gray] (0+\l,0) circle (3pt);
					\draw (4,-0.5) node[inner sep=1pt,rectangle,fill=white] {, };
					\draw (0,-1) node[inner sep=1pt,fill=white] {$\hfill$};
				\end{tikzpicture}
			\end{minipage}
			\\ \hline
	\end{tabular}\end{center}

	\begin{equation}
		\begin{aligned}
			|\mathrm{WHom}^{00}(P_4,P_4 \square P_5)| &=\dbinom{3}{0}|\mathrm{Hom}^{0}(P_1,P_4)| |\mathrm{WHom}^{0}(P_4,P_5)|\\&\quad+\dbinom{3}{1}|\mathrm{Hom}^{0}(P_2,P_4)||\mathrm{WHom}^{0}(P_3,P_5)|\\
			&\quad +\dbinom{3}{2}|\mathrm{Hom}^{0}(P_3,P_4)||\mathrm{WHom}^{0}(P_2,P_5)|\\&\quad+\dbinom{3}{3}|\mathrm{Hom}^{0}(P_4,P_4)|\mathrm{WHom}^{0}(P_1,P_5)|\\
			& =1(1)(13)+3(1)(5)+3(2)(2)+1(3)(1)\\
			&= 43.\nonumber
		\end{aligned}
	\end{equation}

\begin{theorem}\label{Theorem2} Let $m,n$ and $k$ be positive integers and $i,j$ be non-negative integers such that $i < \dfrac{n}{2}-1$ and $j< \dfrac{k}{2}-1$. It follows that
	$$|\mathrm{WHom}^{ij}(P_m,P_n\square P_k)|=\sum_{h=0}^{m-1}\binom{m-1}{h}|\mathrm{Hom}^i(P_{h+1},P_n)||\mathrm{WHom}^j(P_{m-h},P_k)|.$$
\end{theorem}
\begin{proof}
Let $f\in\mathrm{WHom}^{ij}(P_m,P_n\square P_k)$. For each $x\in\{0,1,m-2\}$ in the domain, either $f(x+1)=f(x)\pm(1,0)$ or $f(x+1)=f(x)\pm(0,t)$, where $t\in\{0,1\}$. Assume changes in the first coordinate appear $h$ times. Then, changes in the second coordinate appear $m-1-h$ times. The sequence of changes in the first coordinate form a homomorphism  $f_1\in \mathrm{Hom}^i(P_{h+1},P_n)$. Similarly, the sequence of remaining changes (and no changes) in the second coordinate form a weak homomorphism $f_2\in \mathrm{WHom}^i(P_{m-1-h+1},P_k)$. Thus, the corresponding path graph of $f$ can be obtained from the permutations of all edges in path graphs of $f_1$ and $f_2$ with a fixed sequential order. There are $\binom{m-1}{h}$ permutations in total. Hence, $|\mathrm{WHom}^{ij}(P_m,P_n\square P_k)|=\sum_{h=0}^{m-1}\binom{m-1}{h}|\mathrm{Hom}^i(P_{h+1},P_n)||\mathrm{WHom}^j(P_{m-h},P_k)|.$	
\end{proof}

\noindent From Lemma \ref{Lemma4} and Theorem \ref{Theorem2}, we get the theorem below.\\

\begin{theorem}
	The cardinalities $|\mathrm{WHom}(P_m,P_n\square P_k)|$ of weak homomorphisms from undirected paths $P_m$ to grid graphs $P_n\square P_k$ are	\\
	
	$| \mathrm{WHom}(P_m,P_{n}\square P_{k})|=4\sum_{i=0}^{\lfloor n/2\rfloor-1}\sum_{j=0}^{\lfloor k/2\rfloor-1}| \mathrm{WHom}^{ij}(P_m,P_{n}\square P_{k})| $\\\hspace*{4.3cm}$+(1-(-1)^n)\sum_{j=0}^{\lfloor k/2\rfloor-1}| \mathrm{WHom}^{\lfloor n/2\rfloor j}(P_m,P_{n}\square P_{k})|$\\\hspace*{4.3cm}$+(1-(-1)^k)\sum_{i=0}^{\lfloor n/2\rfloor-1}| \mathrm{WHom}^{i\lfloor k/2\rfloor}(P_m,P_{n}\square P_{k})|$\\\hspace*{3.5cm}$+(1/4)(1-(-1)^n)(1-(-1)^k)|\mathrm{Hom}^{\lfloor n/2\rfloor\lfloor k/2\rfloor}(P_m,P_{n}\square P_{k})|$\\
where\\ 	$|\mathrm{WHom}^{ij}(P_m,P_n\square P_k)|=\sum_{h=0}^{m-1}\binom{m-1}{h}|\mathrm{Hom}^i(P_{h+1},P_n)||\mathrm{WHom}^j(P_{m-h},P_k)|.$
\end{theorem}

For convenience, we compute $|\mathrm{WHom}(P_m,P_n\square P_k)|$ for $2\leq m\leq n,k\leq 8.$ The results are presented in Tables \ref{tab5}.

\begin{table}[H]
	\caption{\label{tab5} Numbers of weak homomorphisms $f:P_m \rightarrow P_n\square P_k$  for $2\leq m\leq n,k\leq 8.$}
	\centering{}%
	\begin{tabular}{c|c|c|c|c|c|c|c|c|}
		\multicolumn{2}{c}{ } & \multicolumn{7}{|c|}{$k$}\\\hline
		$m$&$n$&2&3&4&5&6&7&8\\\hline
		\multirow{7}{*}{2}
		&2&12&20&28&36&44&52&60	\\\cline{2-9}
		&3&20&33&46&59&72&85&98	\\\cline{2-9}
		&4&28&46&64&82&100&118&136	\\\cline{2-9}
		&5&36&59&82&105&128&151&174	\\\cline{2-9}
		&6&44&72&100&128&156&184&212	\\\cline{2-9}
		&7&52&85&118&151&184&217&250	\\\cline{2-9}
		&8&60&98&136&174&212&250&288	\\\hline
		
		\multirow{6}{*}{3}
		&3&&125	&182	& 239 & 296 &353	 &410	\\\cline{2-9}
		&4&&182	&264	&346  &428  &510	 &592	\\\cline{2-9}
		&5&&239	&346	& 453 &560  &667	 &774	\\\cline{2-9}
		&6&&296	&428	& 560 &692  &824	 &956	\\\cline{2-9}
		&7&&353	&510	&667  &824  &981	 &1138	\\\cline{2-9}
		&8&&410	&592	& 774 & 956 &1138	 &1320	\\\hline
		
		\multirow{5}{*}{4} 
		&4&&&1104	&1480	&1856	&2232	&2608	\\\cline{2-9}
		&5&&&1480	&1981	&2482	&2983	&3484	\\\cline{2-9}
		&6&&&1856	&2482	&3108	&3734	&4360	\\\cline{2-9}
		&7&&&2232	&2983	&3734	&4485	&5236	\\\cline{2-9}
		&8&&&2608	&3484	&4360	&5236	&6112	\\\hline
		
		\multirow{4}{*}{5}
		&5&&&&	8733&	11088&13443	&15798	\\\cline{2-9}
		&6&&&&	11088&	14068&	17048&20028	\\\cline{2-9}
		&7&&&&	13443&	17048&	20653&24258	\\\cline{2-9}
		&8&&&&	15798&	20028&	24258&28488	\\\hline
		
		\multirow{3}{*}{6}
		&6&&&&&	 64004&	78226&92448	\\\cline{2-9}
		&7&&&&&	78226 &	95573&112920	\\\cline{2-9}
		&8&&&&&	92448 &112920	&133392	\\\hline
		
		\multirow{2}{*}{7}
		&7&&&&&&	443833&	527452\\\cline{2-9}
		&8&&&&&&	527452&	626696\\\hline

		{8}&8&&&&&&&2951832	\\\hline

	\end{tabular}
	
\end{table}





\end{document}